\documentclass[10pt]{article}
\usepackage{amsmath,amssymb,amsthm,morefloats,todonotes}
\usepackage{epsfig,psfrag,caption}
\usepackage{graphicx,empheq,color,soul}
\usepackage{verbatim}
\usepackage{subfigure}
\numberwithin{equation}{section}

\newcommand{\goto}{\rightarrow}
\newcommand{\bigo}{{\mathcal O}}

\newcommand{\half}{\frac{1}{2}}

\newcommand{\mbf}[1]{\mbox{\boldmath {$#1$}}}

\def\XXint#1#2#3{{\setbox0=\hbox{$#1{#2#3}{\int}$}
     \vcenter{\hbox{$#2#3$}}\kern-.5\wd0}}

\DeclareMathOperator{\diag}{diag}

\DeclareMathOperator{\range}{ran}

\newenvironment{choices}{\left\{ \begin{array}{ll}}{\end{array}\right.}
\newcommand\when{&\text{if~}}
\newcommand\otherwise{&\text{otherwise}}

\newenvironment{mat}{\left[\begin{array}{ccccccccccccccc}}{\end{array}\right]}
\newcommand\bcm{\begin{mat}}
\newcommand\ecm{\end{mat}}

\newcommand{\bea}{\begin{eqnarray}}
\newcommand{\eea}{\end{eqnarray}}
\newcommand{\bean}{\begin{eqnarray*}}
\newcommand{\eean}{\end{eqnarray*}}
\newcommand{\ba}{\begin{array}}
\newcommand{\ea}{\end{array}}
\newcommand{\beqs}{\begin{equation*}\begin{split}}

\newcommand{\mqed}{$\blacksquare$\newline}

\definecolor{DarkGreen}{rgb}{0,.55,0}

\newtheorem{example}{Example}[section]

\newtheorem{remark}{Remark}[section]

\newtheorem{problem}{Problem}[section]

\newtheorem{lemma}{Lemma}[section]

\newtheorem{theorem}{Theorem}[section]

\newtheorem{definition}{Definition}[section]

\newtheorem{assumption}{Assumption}[section]

\newtheorem{corollary}{Corollary}[section]

\newtheorem{proposition}{Proposition}[section]

\newtheorem{algorithm}{Algorithm}[section]

\newcommand{\PII}{Painlev\'e II }

\long\def\symbolfootnote[#1]#2{\begingroup%
\def\thefootnote{\fnsymbol{footnote}}\footnote[#1]{#2}\endgroup}

\def\gmi{i}

\begin{document}
\title{Nonlinear steepest descent and the numerical solution of Riemann--Hilbert problems}
\author{Sheehan Olver$^1$ and Thomas Trogdon$^2$\\
\phantom{.}\\
$^{1}$School of Mathematics and Statistics\\
The University of Sydney\\
NSW 2006, Australia\\
\phantom{.}\\
$^{2}$Department of Applied Mathematics\\
 University of Washington\\
Campus Box 352420\\
 Seattle, WA, 98195, USA \\}
\maketitle

\footnotetext[1]{email: Sheehan.Olver@sydney.edu.au}
\footnotetext[2]{email: trogdon@amath.washington.edu}

\begin{abstract}

The effective and efficient numerical solution of Riemann--Hilbert problems has been
demonstrated in {recent} work.  With the aid of ideas
from the method of nonlinear steepest descent for Riemann--Hilbert
problems, the resulting numerical methods {have been shown numerically to} retain accuracy as values
of certain parameters become arbitrarily large.  Remarkably, this numerical approach does not require knowledge of local parametrices; rather, the deformed contour is scaled near stationary points at a specific rate.   The {primary} aim of this paper
is {to prove that this observed asymptotic accuracy is indeed achieved}.  {To do so, we first construct} a general {theoretical} framework for the
numerical solution of Riemann--Hilbert problems.  Second, we
demonstrate the precise link between nonlinear steepest descent and
the success of numerics in asymptotic regimes.  {In particular, w}e prove sufficient
conditions for numerical methods to retain accuracy.      Finally, we {compute solutions to}
the homogeneous Painlev\'e II equation and the modified Korteweg--de
Vries equations to explicitly demonstrate the {practical validity of the} theory.

\end{abstract}

\section{Introduction}

	Matrix-valued Riemann--Hilbert problems (RHPs) are of profound
        importance in modern applied analysis.    In inverse
        scattering theory,  solutions to the nonlinear Schr\"odinger
        equation, Korteweg de--Vries equation (KdV),
        Kadomtsev--Petviashvili I equation and many other integrable
        solutions can be written in terms of solutions of RHPs
        \cite{AblowitzClarksonSolitons}.   Orthogonal polynomials can also  be rewritten in terms of solutions of RHPs
        \cite{DeiftOrthogonalPolynomials}.  The asymptotics of
        orthogonal polynomials is crucial for determining the
        distribution of eigenvalues of large random matrix ensembles
        \cite{DeiftOrthogonalPolynomials}.   In each of these
        applications, RHPs fulfill the role that integral representations  play in classical asymptotic analysis.

The way {in which} these RHPs are analyzed is through the method of nonlinear steepest descent \cite{DeiftZhouAMS}.  RHPs can be deformed in the complex plane in much the same way as contour integrals.  This allows the oscillatory nature of the problem to be changed to exponential decay.  These deformed problems, which depend on a parameter, are solved approximately.  The approximate solution is found explicitly and the difference between the actual solution and this approximate solution tends to zero as a parameter becomes large.

{The method of nonlinear steepest descent was adapted} by the authors in \cite{TrogdonSOKdV} for
numerical purposes.  We developed a method to reliably solve the
Cauchy initial-value problem for KdV and modified KdV for all values of $x$ and $t$.   {A benefit of the approach was that, unlike the standard method of nonlinear steepest descent, we did not require the knowledge of difficult-to-derive local parametrices.  Instead, an approach was developed based on scaling contours.  Convergence of the approximation was demonstrated through numerical tests, which further showed that the accuracy of the approximation actually increased in the asymptotic regime.}  The focus of the current paper is to derive sufficient conditions (which are often satisfied) for {which we can prove that the approach of solving RHPs numerically on scaled contours}  will be {guaranteed to  be} accurate in asymptotic regimes.  We refer to this type of behavior as asymptotic stability or uniform approximation.

In addition, we show the deep connection between the success of the
numerical method and the success of the method of nonlinear steepest
descent \cite{DeiftZhouAMS}.  A notable conclusion is that one can
expect that  whenever the method of nonlinear steepest descent
produces an asymptotic formula, the numerical method can be made
asymptotically stable.    Achieving this requires varying amounts of
preconditioning of the RHP.  This can vary from not deforming the RHP {at all},
all the way to using the full deformation needed by the analytical
method.  An important question is: ``when can we stop deformations and
a have a reliable numerical method?''  Our main results are in \S \ref{section:uniform} and allow us to answer this question.   {In short, although we do not require the knowledge of local parametrices to construct the numerical method, their existence {ensures} that the numerical method remains accurate, and their explicit knowledge allows us to analyze the error of the approximation directly.}  In the last two sections we provide useful examples of where this arises in applications.

The paper is structured as follows. We begin with some background material, {followed by}  the precise definition of a RHP
  along with properties of an associated singular integral operator.
  This allows us, amongst other things, to address the regularity
  of solutions.
   Then, we use an abstract framework for the numerical
  solution of RHP which will allow us to address asymptotic
  accuracy in a more concise way.  Additionally, other numerical
  methods (besides the one used for the applications) may fit within the
  framework of this paper.  We review the philosophy and analysis behind nonlinear
  steepest descent  and how it relates to our  numerical framework.  Then, we prove our main results which provide sufficient condition for uniform approximation.  The {numerical approach of \cite{SORHFramework} is  placed within the general framework},
  along with  necessary assumptions which allows a realization of uniform
  approximation.  We  apply the theory to two RHPs.  The first is a RHP {representation of the homogeneous} Painlev\'e II transcendent
\begin{align}\label{PII}
u_{xx} - xu-2u^3 = 0,
\end{align}
for specific initial conditions and $x
< 0$.  The second is a RHP {representation of} the modified Korteweg--de Vries equation
\begin{align}\label{mKdV}
u_t - 6u^2u_x + u_{xxx}= 0,
\end{align}
for smooth exponentially decaying initial data in the so-called Painlev\'e
region \cite{deift-zhou:mkdv}.

\section{Background Material}

We use this section to fix notation that will be used throughout the remainder of the manuscript. We reserve $C$ and $C_i$ to denote generic constants, and $X$ and $Y$ for Banach spaces.  We denote the norm on a space $X$ by $\|\cdot\|_X$.  The notation $\mathcal L(X,Y)$ is used to denote the Banach algebra of all bounded linear operators from $X$ to $Y$.  When $X = Y$ we write $\mathcal L(X)$ to simplify notation.  The following lemma is of great use  \cite{atkinson}:
\begin{lemma}\label{op-open} Assume $L \in \mathcal L(X,Y)$ has a bounded inverse $L^{-1} \in \mathcal L(Y,X)$.  Assume $M \in \mathcal L(X,Y)$ satisfies
\begin{align*}
\|M-L\| < \frac{1}{\|L^{-1}\|}.
\end{align*}
Then $M$ is invertible and
\begin{align}\label{op-open-b1}
\|M^{-1}\| \leq \frac{\|L^{-1}\|}{1 - \|L^{-1}\| \|L-M\|}.
\end{align}
Furthermore,
\begin{align}\label{op-open-b2}
\|L^{-1} - M^{-1}\| \leq \frac{\|L^{-1}\|^2 \|L-M\|}{1 - \|L^{-1}\|\|L-M\|}.
\end{align}
\end{lemma}

In what follows we are interested in functions defined on oriented
contours in $\mathbb C$.  Assume $\Gamma$ is piecewise smooth,
oriented and (most likely) self-intersecting. $\gamma_0$ is used to
denote the set of self-intersections.  Decompose $\Gamma = \Gamma_1
\cup \cdots \cup \Gamma_{\ell}$ to its smooth, non-self-intersecting
components.  Define $C^\infty_c(\Gamma_\gmi)$ to be the space of
infinitely differentiable functions with compact support in $\Gamma_\gmi$, and $C(\Gamma_\gmi)$ to be the Banach space of continuous functions with the uniform norm.  Define the space
\begin{align*}
L^2(\Gamma) = \left\{f \mbox{ measurable}:  \sum_{\gmi=1}^\ell \int_{\Gamma_\gmi} |f(k)|^2 |dk| < \infty \right\}.
\end{align*}

We define, $D$, the distributional differentiation operator for functions defined on $\Gamma\setminus \gamma_0$.  For a function $\varphi \in C^\infty_c(\Gamma_\gmi)$ we represent a linear function $g$ via the dual pairing  
\begin{align*}
g(\varphi) = \langle g, \varphi \rangle_{\Gamma_\gmi}.
\end{align*}
 To $D_\gmi g$ we associate the functional
\begin{align*}
\langle g, \varphi'_\gmi \rangle_{\Gamma_\gmi}.
\end{align*}
For $f \in L^2(\Gamma)$ consider $f_\gmi = f|_{\Gamma_\gmi}$, the restriction of $f$ to $\Gamma_\gmi$.  In the case that the distribution $D_\gmi f_\gmi $ corresponds to a locally integrable function  
\begin{align*}
\langle D_\gmi f_\gmi, \varphi \rangle_{\Gamma_\gmi} = \int_{\Gamma_\gmi} D_\gmi f_\gmi(k) \varphi(k) dk = \int_{\Gamma_\gmi} f_\gmi(k) \varphi'(k) dk,
\end{align*}
for each $i$, we define
\begin{align*}
Df(k) = D_\gmi f_\gmi(k) \mbox{ if } k \in \Gamma_\gmi\setminus \gamma_0.
\end{align*}
This allows us to define
\begin{align*}
H^k(\Gamma) = \left\{ f \in L^2(\Gamma) : D^jf \in L^2(\Gamma), ~~j=0,\ldots,k\right\}.
\end{align*}
We write $W^{k,\infty}(\Gamma)$ for the Sobolev space  with the $L^2$ norm replaced with the $L^\infty$ norm.  An important note is that we will be dealing with matrix-valued functions, and hence the definitions of all these spaces must be suitably extended.  Since all finite-dimensional norms are equivalent, we can use the above definitions in conjunction with any matrix norm to define a norm for matrix-valued functions provided the norm is sub-additive.

\section{Theory of Riemann--Hilbert Problems in $L^2(\Gamma)$} \label{section:L2}

Loosely speaking, a Riemann--Hilbert problem (RHP) is a boundary-value problem in the complex plane:

\begin{problem}\label{RHproblem} \cite{SIE} Given an oriented contour $\Gamma \subset
  \mathbb C$ and a jump matrix $G: \Gamma \goto \mathbb C^{2\times
    2}$, find a bounded function $\Phi: \mathbb C \setminus \Gamma
  \goto \mathbb C^{2\times 2}$ which is analytic everywhere in the
  complex plane except on $\Gamma$, such that 
\begin{align}
\Phi^+(z) &= \Phi^-(z) G(z), ~~~ \mbox{ for $z \in \Gamma$, and } \label{bvp}\\
\Phi(\infty) &= I, \label{infty}
\end{align}
where $\Phi^+$ denotes the limit of $\Phi$ as $z$ approaches $\Gamma$
from the left, $\Phi^-$ denotes the limit of $\Phi$ as $z$
approaches $\Gamma$ from the right, and $\Phi(\infty) = \lim_{|z| \goto
  \infty} \Phi(z)$.   We denote this RHP by $[G;\Gamma]$.
\end{problem}

The definition is not
sufficiently precise to compute solutions.  We say
that $\Phi$ is an $L^2$ solution of $[G;\Gamma]$ normalized at
infinity, if
\begin{align*}
\Phi(z) = I + \frac{1}{2\pi i} \int_{\Gamma} \frac{ u(s) }{s-z} ds = I + \mathcal C_\Gamma u(z),
\end{align*}
 for some $u \in L^2(\Gamma)$  and $\Phi^+(z) = \Phi^-(z) G(z)$, for
 almost every $z \in \Gamma$.

\begin{definition}
The \emph{Cauchy transform pair} is the pair of boundary values
of the Cauchy integral $\mathcal C_\Gamma$:
\begin{align}\label{cauchy-trans}
\mathcal C_\Gamma^\pm u(z) = \left(\mathcal C_\Gamma u(z) \right)^\pm .
\end{align}
\end{definition}

It is well known that for fairly general contours (including every contour we
considered here) that the limits in \eqref{cauchy-trans} exist
non-tangentially almost everywhere.  The operators
$\mathcal C_\Gamma^\pm$ are
bounded from $L^2(\Gamma)$ to itself and satisfy the operator identity \cite{wavelets}
\begin{align}\label{proj}
\mathcal C_\Gamma^+ - \mathcal C_\Gamma^- = I.
\end{align}

\begin{remark}  When deriving a RHP one must show that $\Phi - I \in
  \range \mathcal C_\Gamma$.  This is generally done using Hardy space
  methods \cite{wavelets,Duren}.  Brevity is the governing motivation for the above definition of a RHP.
\end{remark}

We convert the RHP into an equivalent singular integral
equation (SIE).  Assume $\Phi(z) = I + \mathcal C_\Gamma u(z)$ and
substitute into \eqref{bvp} to obtain
\begin{align}
I + \mathcal C_\Gamma^+u(z) = \mathcal C_\Gamma^- u(z) G(z) + G(z).
\end{align}
Using \eqref{proj},
\begin{align}\label{SIE}
u(z) - \mathcal C_\Gamma^-u(z) (G(z) - I) = G(z) - I.
\end{align}
\begin{definition}
 	We use $\mathcal C [G;\Gamma]$ to refer to the operator
(bounded on $L^2(\Gamma)$ provided $G \in L^\infty(\Gamma)$) in
\eqref{SIE}.
\end{definition}

  In what follows we assume at a minimum that the RHP
is well-posed, or $\mathcal C[G;\Gamma]^{-1}$ exists and is bounded in
$L^2(\Gamma)$ and $G-I \in L^2(\Gamma)$.

We need to establish the smoothness of solutions of \eqref{SIE} since
we  approximate solutions numerically.  This smoothness relies
on the smoothness of the jump matrix in the classical sense along with a type of smoothness
at the self-intersection points, $\gamma_0$, of the contour.

\begin{definition}\label{prod-cond}  Assume $a \in \gamma_0$. Let $  \Gamma_1, \ldots,   \Gamma_m$ be a
counter-clockwise ordering of subcomponents of $\Gamma$ which contain $z = a$ as
an endpoint.  For $G \in W^{k,\infty}(\Gamma)\cap H^k(\Gamma)$ we define $\hat G_i$ by $G|_{\tilde \Gamma_i}$ if $ \Gamma_i$ is oriented outwards and $(G|_{ \Gamma_i})^{-1}$ otherwise.  We say $G$ satisfies the $(k-1)$th-order product condition if using the $(k-1)$th-order Taylor expansion {of each $\hat G_i$} we have
\begin{align}
\prod_{i=1}^m \hat G_i = I + \bigo\left((z-a)^k\right) , \mbox{ for } j = 1, \ldots, k-1, ~~ \forall a \in \gamma_0.
\end{align}
\end{definition}

\begin{figure}[ht]
\begin{center}
\includegraphics[width=.5\linewidth]{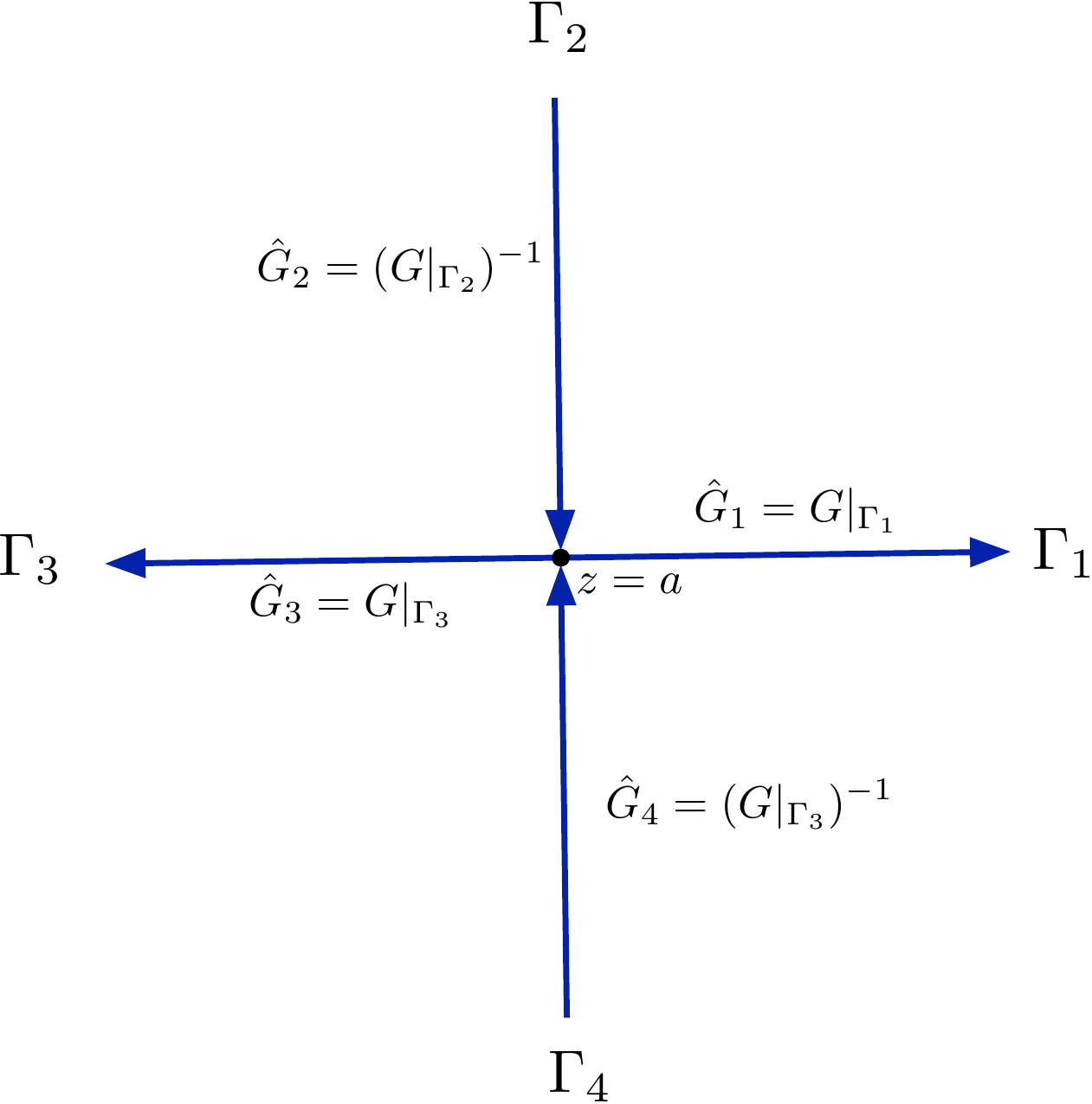}
\caption{\label{productcond} A typical intersection point and the
  definition of $\hat G_i$.}
\end{center}
\end{figure}

To capture all the regularity properties of the solution of
\eqref{SIE}, and to aid the development of numerical methods, we
introduce the following definition in the same vein as Definition \ref{prod-cond}.

\begin{definition}\label{zero-sum} Assume that $a \in \gamma_0$ and let $ \Gamma_1, \ldots,  \Gamma_m$ be a
counter-clockwise ordering of subcomponents of $\Gamma$ which contain $z = a$ as
an endpoint.  For $f \in H^{k}(\Gamma)$, define
\begin{align}
f^{(j)}_i = \begin{choices} -\lim_{z \goto a} \left(\frac{d}{dz}\right)^j f|_{ \Gamma_i}(z)
  \when  \Gamma_i \mbox{ is oriented outward,}\\
\lim_{z \goto a} \left(\frac{d}{dz}\right)^j f|_{ \Gamma_i}(z)
  \when  \Gamma_i \mbox{ is oriented inward.} \end{choices}
\end{align}
We say that $f$ satisfies the $(k-1)$th-order zero-sum condition if
\begin{align}
\sum_{i=1}^m f^{(j)}_i = 0, \mbox{ for } j = 0, \ldots, k-1 \mbox{ and
} \forall a \in \gamma_0.
\end{align}
\end{definition}
\begin{remark}  These definitions imply conditions also when $\Gamma$ has an isolated endpoint. In that case $G = I$ at this endpoint and $f = 0$.
\end{remark}

This motivates the following definition.
\begin{definition}  A RHP $[G;\Gamma]$ is said to be $k$-regular if 
\begin{itemize}
\item $G, G^{-1}$ satisfy the $(k-1)$th-order product condition, and
\item $G- I, G^{-1} -I \in W^{k,\infty}(\Gamma) \cap H^k(\Gamma)$.
\end{itemize}
\end{definition}
It is worth noting that when $\Gamma$ is bounded the $H^k(\Gamma)$ condition is trivial.  We include it for completeness.  A useful result \cite[Chapter II, Lemma 6.1]{mikhlin} is:

\begin{lemma}\label{mikhlin} Let $\Gamma$ be a smooth non-closed curve from $z =a$ to $z =b$ with orientation from $a$ to $b$.  If $f \in H^{1}(\Gamma)$, then
\begin{align*}
D \int_{\Gamma} \frac{f(\zeta)}{\zeta-z} d\zeta = \frac{f(a)}{a-z} - \frac{f(b)}{b-z} + \int_{\Gamma} \frac{Df(\zeta)}{\zeta-z}d\zeta.
\end{align*}
\end{lemma}

{Due to cancellation}, $\mathcal
C_\Gamma^\pm$ commutes with the weak differentiation operator for functions that satisfy the zeroth-order zero-sum condition.  We use the notation $H_z^{k}(\Gamma)$ for the closed subspace
of $H^{k}(\Gamma)$ consisting of functions which satisfy the $(k-1)$th
order zero-sum condition. This allows us to use the boundedness of $\mathcal C_\Gamma^\pm$ on $L^2(\Gamma)$ to show the boundedness (with same norm) of $\mathcal C_\Gamma^\pm$ from $H^k_z(\Gamma)$ to $H^k(\Gamma)$.  The results in \cite{zhou-RHP}, combined
with the fact that differentiation commutes with the Cauchy operator on these zero-sum spaces 
can be used to prove the following theorem.

\begin{theorem} \label{regularity}
Given a RHP $[G;\Gamma]$ which is $k$-regular, assume  $\mathcal C[G;\Gamma]$ is invertible on $L^2(\Gamma)$. Then $u \in H^k_z(\Gamma)$, the solution of \eqref{SIE}, satisfies
\begin{align}\label{uk}
D^ku = \mathcal C[G;\Gamma]^{-1} \biggl( D^k(G-I) + D^k( \mathcal C_\Gamma^- u \cdot (G-I)) -  \mathcal C_\Gamma^{-} D^ku \cdot (G-I)\biggr).
\end{align}
where the right-hand side of \eqref{uk} does not depend on $D^ku$.
\end{theorem}

\begin{corollary}  \label{regularity-cor} Under the hypotheses of Theorem \ref{regularity},
  $u$ satisfies an inequality of the form
\begin{align}\label{hk-bound}
\|u\|_{H^k(\Gamma)} \leq p_k\left(\|G-I\|_{W^{k,\infty}(\Gamma)}\|\mathcal C[G;\Gamma]^{-1}\|_{\mathcal L(L^2(\Gamma))}\right) \|\mathcal C[G;\Gamma]^{-1}\|_{\mathcal L(L^2(\Gamma))}\|G-I\|_{H^k(\Gamma)},
\end{align}
where $p_k$ is a polynomial of degree $k$ whose coefficients depend on  $\|\mathcal C_\Gamma^-\|_{\mathcal L(L^2(\Gamma))}$.
\end{corollary}
\noindent\emph{Proof:} Taking the norm of \eqref{uk} gives a bound on the semi-norm $\|D^ku\|_{L^2(\Gamma)}$ in terms of $\|u\|_{H^{k-1}(\Gamma)}$.  Using \eqref{uk} for $k-1$ gives a bound in terms of $\|u\|_{H^{k-2}(\Gamma)}$.  This process produces a bound of the form \eqref{hk-bound}. \mqed

\begin{remark} The expression for the derivative in Theorem \ref{regularity} can be used to bound Sobolev norms of the solution in terms of Sobolev norms of the jump matrix if a bound on the norm of the inverse operator is known.  In some cases, when the jump matrix depends on a parameter, and the Sobolev norms of the jump matrix are bounded or decaying, the resulting bounds are of use.  
\end{remark}

 We often use the following theorem
which is derived from the results in \cite{FokasPainleve}.

\begin{theorem} \label{approx-prop}
Consider a sequence of RHPs $\{[G_\xi;\Gamma]\}_{\xi \geq 0}$ on
the fixed contour $\Gamma$ which are $k$-regular.  Assume $G_\xi \goto G$ in
$W^{k,\infty}(\Gamma) \cap H^k(\Gamma)$ as $\xi \goto \infty$, then
\begin{itemize}
\item If $\mathcal C[G;\Gamma]$ is invertible then there exists a $T>0$
such that $\mathcal C [G_\xi;\Gamma]$
is also invertible for $\xi > T$.
\item If $\Phi_\xi$ is the solution of $[G_\xi,\Gamma]$, and $\Phi$ the
  solution of $[G;\Gamma]$, then $\Phi_\xi^\pm - \Phi^\pm \goto 0$ in
  $H^k(\Gamma)$.
\item $\| \Phi_\xi - \Phi\|_{W^{j,\infty}(S)} \goto 0$, for all $j \geq 1$, whenever $S$ is bounded away from $\Gamma$.
\end{itemize}
\end{theorem}

\noindent \emph{Proof: }  The first statement follows from the fact that $\mathcal C[G_\xi;\Gamma]$ converges to $\mathcal C[G;\Gamma]$ in operator norm. The second property follows from Corollary \ref{regularity-cor}.  The final property is a consequence of the Cauchy-Schwartz inequality and the fact that $\|u_\xi-u\|_{L^2(\Gamma)} \goto 0$. \mqed

\section{The Numerical Solution of Riemann--Hilbert Problems}\label{section:setting}

The goal in this section is to introduce the necessary tools to approximate the operator equation
\begin{align}\label{SIE-abstract}
\mathcal C[G;\Gamma]u = G-I.
\end{align}
We start with two projections $\mathcal I_n$ and $\mathcal P_n$, both of finite rank.  Assume $\mathcal P_n$ is defined on $H^1_z(\Gamma)$ and $\mathcal I_n$ is defined on $H^1(\Gamma)$.  Define $X_n = \range \mathcal P_n$ and $Y_n = \range \mathcal I_n$ equipping both spaces with the inherited $L^2(\Gamma)$ norm.  We obtain a finite-dimensional approximation of $\mathcal C[G;\Gamma]$ by defining
\begin{align*}
\mathcal C_n[G;\Gamma]u = \mathcal I_n \mathcal C[G;\Gamma] u.
\end{align*}
It follows that $\mathcal C_n[G;\Gamma]: X_n \goto Y_n$.  This is useful under the assumption that we can compute $\mathcal C[G;\Gamma]$ exactly for all $u \in X_n$.  An approximate solution of \eqref{SIE-abstract} is obtained by
\begin{align*}
u_n = \mathcal C_n[G;\Gamma]^{-1} \mathcal I_n (G-I),
\end{align*}
whenever the operator is invertible.  We use the pair $(\mathcal I_n,\mathcal P_n)$ to refer to this numerical method.

\begin{remark}
 One may have concerns about $\mathcal C_n[G;\Gamma]$ because the dimension of the domain and that of the range are different.  It turns out that $\mathcal C[G;\Gamma]$ maps $H_z^1(\Gamma)$ to a closed subspace of $H^1(\Gamma)$ and we can define $\mathcal I_n$ on this space.  In {the numerical framework of \cite{SORHFramework}, solving the associated linear system  results in a solution that must satisfy the zeroth-order zero-sum condition, justifying the theoretical construction above}. In what follows we ignore this detail.
\end{remark}

To simplify notation, we define $\mathcal T[G;\Gamma]u =\mathcal C_\Gamma^- u (G-I)$ {(so that $\mathcal C[G;\Gamma] = I - \mathcal T[G;\Gamma]$)}  and $\mathcal T_n[G;\Gamma] = \mathcal I_n \mathcal T[G;\Gamma]$.  We use a few definitions to describe the required properties of the projections.

\newcommand{\Norms}{\|G-I\|_{L^\infty(\Gamma)} \|\mathcal C_\Gamma^-\|_{\mathcal L(L^2(\Gamma))}}

\begin{definition}\label{alpha-beta}  The approximation $\mathcal C_n[G;\Gamma]$ to $\mathcal C[G;\Gamma]$ is said to be of type $(\alpha,\beta,\gamma)$ if, whenever $\mathcal C[G;\Gamma]$ is invertible for $n > N$, $\mathcal C_n[G;\Gamma]$ is invertible and
\begin{itemize}
\item $\|\mathcal C_n[G;\Gamma]\|_{\mathcal L(H_z^1(\Gamma),Y_n)} \leq C_1 n^\alpha (1 + \Norms)$,
\item $\|\mathcal C_n[G;\Gamma]^{-1}\|_{\mathcal L(Y_n,X_n)} \leq C_2
  n^\beta \|\mathcal C[G;\Gamma]^{-1}\|_{\mathcal L(L^2(\Gamma))}$ and
\item $\|\mathcal T_n[G;\Gamma]\|_{\mathcal L(X_n,Y_n)} \leq C_3 n^\gamma \Norms$.
\end{itemize}
The constants here are allowed to depend on $\Gamma$.
\end{definition}

The first and second conditions in Definition \ref{alpha-beta} are necessary for the convergence of the numerical method.  This will be made more precise below.  The first and third conditions are needed to control operator norms as $G$ changes.  It is not surprising that the first and the third conditions are intimately related and in \S \ref{section:realize} we demonstrate the connection.

\begin{remark} Some projections, mainly those used in Galerkin methods, can be defined directly on $L^2(\Gamma)$.  In this case we replace the first condition in Definition \ref{alpha-beta} with 
\begin{align*}
 \|\mathcal C_n[G;\Gamma]\|_{\mathcal L(L^2(\Gamma),Y_n)} \leq C_1 n^\alpha (1+\Norms).
\end{align*}
This condition and the second condition with $\alpha = \gamma$ are implied by requiring
\begin{align*}
\|\mathcal I_n\|_{\mathcal L(L^2(\Gamma),Y_n)} \leq C_1 n^\alpha.
\end{align*}
In this sense Galerkin methods are more natural for RHPs, though we use {the  collocation method of \cite{SORHFramework} below because a Galerkin method has yet to be developed.}
\end{remark}

\begin{definition}\label{admissible}
The pair $(\mathcal I_n, \mathcal P_n)$ is said
to produce an admissible numerical method if 
\begin{itemize}
\item The method is of type $(\alpha,\beta,\gamma)$.
\item For all $m> 0$, $\|\mathcal I_n u -u \|_{H^{1}(\Gamma)}$ and $\|\mathcal P_n u -u
  \|_{H^{1}(\Gamma)}$ tend to zero faster than $n^{-m}$ as $n \goto \infty$ for all $u \in H^s(\Gamma)$; for some $s$.
\item $\mathcal I_n$ is bounded from $C(\Gamma)$ to $L^2(\Gamma)$, uniformly in $n$.
\end{itemize}
\end{definition}

\begin{remark}  We  assume spectral convergence of the projections.
  This can be relaxed but one has to spend considerable effort to
  ensure $\alpha$, $\beta$ and $\gamma$ are sufficiently small.
\end{remark}

Next, we prove the generalized convergence theorem.

\begin{theorem}\label{convergence}
Assume that $(\mathcal I_n, \mathcal P_n)$ produces an admissible
numerical method.  If $\mathcal C[G;\Gamma]$ is 1-regular and invertible on
$L^2(\Gamma)$, we have
\begin{align}\label{L2-estimate}
\|u-u_n\|_{L^2(\Gamma)} &\leq (1 + cn^{\alpha+\beta}) \|\mathcal P_n u -
u\|_{H^1(\Gamma)}~~ \mbox{ with }\\
c &= C \|\mathcal C [G;\Gamma]^{-1}\|_{\mathcal L(L^2(\Gamma))}(1+\Norms).
\end{align}
\end{theorem}

\noindent \emph{Proof: } First, for notational simplicity, define $\mathcal K_n
= \mathcal C_n[G;\Gamma]$, $\mathcal K = \mathcal C[G;\Gamma]$ and $f
= G-I$.  Then $u_n = \mathcal K_n^{-1}
\mathcal I_n f = \mathcal K_n^{-1} \mathcal I_n \mathcal K u$.  Further, since $u \in H^1_z(\Gamma)$,
\begin{align*}
u - u_n =& u - \mathcal P_n u + \mathcal P_n u - u_n\\
=& u - \mathcal P_n u + \mathcal P_n u - \mathcal K_n^{-1}\mathcal I_n
\mathcal K u\\
=& u - \mathcal P_n u+ \mathcal K_n^{-1} \mathcal K_n \mathcal P_n u-
\mathcal K_n^{-1} \mathcal
I_n \mathcal K u\\
=& u-\mathcal P_n u+ \mathcal K_n^{-1}(\mathcal K_n \mathcal P_n u - \mathcal I_n \mathcal K u)\\
=& u-\mathcal P_n u + \mathcal K_n^{-1} \mathcal I_n \mathcal K(\mathcal P_n u - u).
\end{align*}
We used $\mathcal K_n \mathcal P_n u = \mathcal I_n \mathcal K \mathcal P_n u$ for  $u \in H^1_z(\Gamma)$ in the last line.  Taking an $L^2(\Gamma)$ norm, we have
\begin{align*}
\|u-u_n\|_{L^2(\Gamma)} \leq \|(I + \mathcal K_n^{-1} \mathcal I_n
\mathcal K(u-\mathcal P_nu) \|_{L^2(\Gamma)} \leq (1+cn^{\alpha+\beta})
\|u-\mathcal P_n u\|_{H^1(\Gamma)}. \blacksquare
\end{align*}

\begin{remark}  In the case mentioned above where $\mathcal I_n$ and $\mathcal P_n$ can be defined directly on $L^2(\Gamma)$ we obtain a purely $L^2(\Gamma)$ based bound
\begin{align*}
\|u-u_n\|_{L^2(\Gamma)} \leq (1+cn^{\alpha+\beta})
\|u-\mathcal P_n u\|_{L^2(\Gamma)}.
\end{align*}
\end{remark}

\begin{corollary}\label{convergence-rhp}
Under the assumptions of Theorem \ref{convergence} and assuming that
$[G;\Gamma]$ is $k$-regular for large $k$ (large is determined by Definition \ref{admissible}) we have that
$\Phi_n = I + \mathcal C_\Gamma u_n$ is an approximation of $\Phi$, the solution of $[G;\Gamma]$,
in the following sense.
\begin{itemize}
\item $\Phi_n^\pm - \Phi^\pm \goto 0 $ in $L^2(\Gamma)$ and
\item $\|\Phi_n - \Phi\|_{W^{j,\infty}(S)} \goto 0$ for all $j \geq 0$,  whenever $S$ is bounded away from  $\Gamma$.
\end{itemize}
\end{corollary}

\noindent \emph{Proof: } The first claim follows from the boundedness
of the Cauchy operator on $L^2(\Gamma)$ and, as before, the
Cauchy-Schwartz inequality gives the second. \mqed

Below we always assume the numerical method considered is admissible.  The ideas presented thus far are general.  In specific cases the
contour $\Gamma$ consists of disjoint components.  We take a different
approach to solving the RHP in this case.  

\begin{example}  \label{example:disjoint} Consider the RHP $[G;\Gamma]$ with $\Gamma= \Gamma_1
  \cup \Gamma_2$ where $\Gamma_1$ and $\Gamma_2$ are disjoint.  To
  {solve the full RHP, we first  solve  for $\Phi_1$ --- the solution of
  $[G|_{\Gamma_1};\Gamma_1]$ --- assuming that this sub-problem has a unique solution}.  The jump on $\Gamma_2$ is modified
   through conjugation by $\Phi_1$. Define
\begin{align*}
\tilde G_2 = \Phi_1 G|_{\Gamma_2} \Phi_1^{-1}.
\end{align*}
The solution $\Phi_2$ of $[\tilde G_2;\Gamma_2]$ is then found.  A
simple calculation shows that $\Phi = \Phi_1\Phi_2$ {solves the original RHP}
$[G;\Gamma]$. \end{example}

This idea allows us to treat each disjoint contour separately, solving in an iterative way.
When using this algorithm numerically, the dimension of the linear
system solved at each step is a fraction of that of the full discretized problem.  This
produces significant computational savings.  We now generalize these ideas.

Consider a RHP $[G;\Gamma]$ where $\Gamma = \Gamma_1 \cup
\cdots \cup \Gamma_\ell$.  Here each $\Gamma_\gmi$ is disjoint and $\Gamma_\gmi
= \alpha_\gmi \Omega_\gmi + \beta_\gmi$ for some contour $\Omega_\gmi$.  We define
$G_\gmi(z) = G(z)|_{\Gamma_\gmi}$ and $H_\gmi(k) = G_\gmi(\alpha_\gmi k + \beta_\gmi)$.  As a notational remark, we always associate $H_i$ and $G$ is this way.

\begin{remark}
The motivation for introducing {the representation of the contours} in this fashion will be made
clear below.  Mainly, this formulation is important when $\alpha_i$
and/or $\beta_i$ depend on a parameter but $\Omega_i$ does not.
\end{remark}

We now describe the general iterative solver.
\begin{algorithm} \label{scale-shift} (Scaled and Shifted RH Solver)
\begin{enumerate}
\item Solve the RHP $[H_1; \Omega_1]$ to obtain $\tilde
  \Phi_1$.  We denote the solution of the associated SIE as $U_1$ with
  domain $\Omega_1$.  Define $\Phi_1(z) = \tilde \Phi_1
  \left(\frac{z-\beta_1}{\alpha_1} \right)$.
\item For each $j = 2, \ldots, \ell$ define $\Phi_{i,j}(z) =
  \Phi_i(\alpha_j z + \beta_j)$ and solve the RHP
$[\tilde H_j; \Omega_j]$ with 
\begin{align*}
\tilde H_j = \Phi_{j-1,j} \cdots \Phi_{1,j} H_j \Phi_{1,j}^{-1} \cdots
  \Phi_{j-1,j}^{-1},
\end{align*}
 to obtain $\tilde \Phi_j$.  Again, the
  solution of the integral equation is denoted by $U_j$ with domain
  $\Omega_j$.  Define $\Phi_j(z) = \tilde
  \Phi_j\left(\frac{z-\beta_j}{\alpha_j}\right)$.
\item Construct $\Phi = \Phi_\ell \cdots \Phi_1$, which satisfies the
  original problem.
\end{enumerate}
\end{algorithm}

When this algorithm is implemented numerically, the jump matrix corresponding to $\tilde H_j$ is not exact.  It  depends on the approximations of each of the $\Phi_i$ for $i< j$ and more specifically, it  depends on the order of approximation of the RHP on $\Omega_i$ for $i<j$. We use the notation $\mbf n_i = (n_1, \ldots, n_i)$ where each $n_i$ is the order of approximation on $\Omega_i$. We use $\mbf n > \mbf m$ whenever
the vectors are of the same length and $n_j > m_j$ for all $j$.  The
statement $\mbf n \goto \infty$ means that each component of $\mbf n$ tends to $\infty$.  Let $\Phi_{i,j,\mbf n_i}$ be the approximation of $\Phi_{i,j}$ and define
\begin{align*}
\tilde H_{j,\mbf n_j} = \Phi_{j-1,j,\mbf n_{j-1}} \cdots \Phi_{1,j,\mbf n_1} H_j \Phi_{1,j,\mbf n_1}^{-1} \cdots
  \Phi_{j-1,j,\mbf n_{j-1}}^{-1}.
\end{align*}
If the method converges then $\tilde H_{j,\mbf n_j} \goto \tilde H_j$ uniformly as $\mbf n_j \goto \infty$.

A significant remaining question is: ``how do we know solutions exist at
each stage of this algorithm?'' In general, this is not the case.  $\mathcal C[G;\Gamma]$ can be expressed in the form
$\mathcal K - \mathcal T$ where $\mathcal K$ is the block-diagonal
operator with blocks $\mathcal C[G_i;\Gamma_i]$ and $\mathcal T$ is a
compact operator.  Here $\mathcal T$ represents the effect of one
contour on another and if the operator norm of $\mathcal T$ is
sufficiently small
solutions exist at each iteration of Algorithm
\ref{scale-shift}.  This is true if the arclength of each
$\Gamma_i$ is sufficiently small.  We leave a more thorough discussion of this to \S
\ref{subsection:nsd}.  An implicit assumption in our numerical {framework} is that
such equations are uniquely solvable.

The final question is one of convergence.  For a single fixed contour we know that if $(\mathcal I_n, \mathcal P_n)$ produces an admissible numerical method and the RHP is sufficiently regular, the numerical method  converges.  This means that the solution of this RHP  converges uniformly, away from the contour it is defined on.  This is the basis for proving that Algorithm \ref{scale-shift} converges. Theorem \ref{approx-prop} aids us when considering the infinite-dimensional operator {for which} the jump matrix is uniformly close, but we need an {additional} result for the finite-dimensional case.

\begin{lemma} \label{discrete-approx-prop}  Consider a sequence of RHPs $\{[G_\xi;\Gamma]\}_{\xi \geq 0}$ on
the fixed contour $\Gamma$ which are $k$-regular.  Assume $G_\xi \goto G$ in
$L^{\infty}(\Gamma) \cap L^2(\Gamma)$ as $\xi \goto \infty$ and $[G;\Gamma]$ is $k$-regular, then
\begin{itemize}
\item If $\mathcal C_n[G;\Gamma]$ is invertible, then there exists $T(n)>0$
such that $\mathcal C_n [G_\xi;\Gamma]$ is also invertible for $\xi > T(n)$.
\item If $\Phi_{n,\xi}$ is the approximate solution of $[G_\xi;\Gamma]$ and $\Phi_n$ is the approximate
  solution of $[G;\Gamma]$, then $\Phi_{n,\xi} - \Phi_n \goto 0$ in
  $L^2(\Gamma)$ as $\xi \goto \infty$ for fixed $n$.
\item $\| \Phi_{n,\xi} - \Phi_n\|_{W^{j,\infty}(S)} \goto 0$, as $\xi \goto \infty$, for all $j \geq 1$, whenever $S$ is bounded away from $\Gamma$ for fixed $n$.
\end{itemize}
\end{lemma}

\noindent \emph{Proof: }  We consider the two equations
\begin{align*}
\mathcal C_n[G_\xi;\Gamma] u_{n,\xi} = \mathcal I_n(G_\xi-I),\\
\mathcal C_n[G;\Gamma] u_n = \mathcal I_n (G-I).
\end{align*}
Since the method is of type $(\alpha,\beta,\gamma)$, we have (see Definition \eqref{alpha-beta}),
\begin{align*}
\|\mathcal C_n[G_\xi;\Gamma] - \mathcal C_n[G;\Gamma]\|_{\mathcal L(X_n,Y_n)} \leq C_3 n^\gamma \|\mathcal C_\Gamma^-\|_{\mathcal L(L^2(\Gamma))} \|G_\xi-G\|_{L^\infty(\Gamma)} = E(\xi) n^{\gamma}.
\end{align*}
For fixed $n$, by increasing $\xi$, we can make $E(\xi)$ small, so that
\begin{align*}
 \|\mathcal C_n[G_\xi;\Gamma] - \mathcal C_n[G;\Gamma]\|_{\mathcal L(X_n,Y_n)}\leq \frac{1}{C_2}\frac{1}{ \|\mathcal C[G;\Gamma]^{-1}\|_{\mathcal L(L^2(\Gamma))}} n^{-\beta} \leq \frac{1}{\|\mathcal C_n[G;\Gamma]^{-1}\|_{\mathcal L(Y_n,X_n)}}.
\end{align*}
Specifically, we choose $\xi$ small enough so that
\begin{align*}
E(\xi) \leq \half \frac{1}{C_2C_3}\frac{1}{ \|\mathcal C[G;\Gamma]^{-1}\|_{\mathcal L(L^2(\Gamma))}} n^{-\gamma-\beta}.
\end{align*}
Using Lemma \ref{op-open} $\mathcal C_n[G_\xi;\Gamma]$ is invertible, and we bound
\begin{align}\label{discrete-bound}
\|\mathcal C_n[G_\xi;\Gamma]^{-1} - \mathcal C_n[G;\Gamma]^{-1}\|_{\mathcal L(Y_n,X_n)} &\leq 2C_2 n^{2\beta+\gamma} \|\mathcal C[G;\Gamma]^{-1}\|_{\mathcal L(L^2(\Gamma))}^2 E(\xi).
\end{align}
Importantly, the quantity on the left tends to zero as $\xi \goto \infty$.  We use a triangle inequality
\begin{align*}
\|u_n-u_{n,\xi}\|_{L^2(\Gamma)} \leq \|(\mathcal C_n[G_\xi;\Gamma]^{-1} - \mathcal C_n[G;\Gamma]^{-1})\mathcal I_n (G-I)\|_{L^2(\Gamma)} + \| \mathcal C_n[G;\Gamma]^{-1}\mathcal I_n (G-G_\xi)\|_{L^2(\Gamma)}.
\end{align*}
Since we have assumed that $\Gamma$ is bounded and that the norm of
$\mathcal I_n: C(\Gamma) \goto L^2(\Gamma)$ is uniformly bounded in $n$, we obtain $L^2$ convergence of $u_n$ to $u_{n,\xi}$ as $\xi \goto \infty$:
\begin{align}\label{discrete-final-bound}
\|u_n-u_{n,\xi}\|_{L^2(\Gamma)} \leq C_3n^{2\beta+\gamma}E(\xi)\|G-I\|_{L^\infty(\Gamma)} + C_4 n^\beta \|G-G_\xi\|_{L^\infty(\Gamma)} \leq C_5n^{2\beta+\gamma} E(\xi).
\end{align}
This proves the three required properties. \mqed

\begin{remark} A good way to interpret this result is to see $E(\xi)$
  as the difference in norm between the associated infinite-dimensional operator which is proportional to the uniform difference in the
  jump matrices. Then \eqref{discrete-bound} gives the resulting error
  between the finite-dimensional operators. 
  It is worthwhile to note that if $\alpha = \beta = \gamma = 0$ then $\delta$ can be chosen independent of $n$.
\end{remark}

Now we have the tools needed to address the convergence of the solver.
We introduce some notation to simplify {matters}.  At stage $j$ in the
solver we solve a SIE on $\Omega_j$.  On this domain we need to compare two RHPs:
\begin{itemize}
\item $[\tilde H_j;\Omega_j]$ and
\item $[\tilde H_{j,\mbf n_j};\Omega_j]$.
\end{itemize}
Let $U_j$ be the exact solution
of this SIE which is obtained from $[\tilde H_j;\Omega_j]$.  As an intermediate step we need to consider the numerical solution of $[\tilde H_j;\Omega_j]$.  We use $U_{j,n_j}$ to denote the numerical approximation of $U_j$ of order $n_j$.  Also, $U_{j,\mbf n_j}$ will be used to denote the numerical approximation of the solution of the SIE associated with $[\tilde
H_{j,\mbf n_j};\Omega_j]$.

\begin{theorem}  \label{algorithm-convergence} Assume that each problem in Algorithm \ref{scale-shift} is solvable and $k$-regular for sufficiently large $k$. Then the algorithm converges to the true solution of the RHP.  More precisely, there exists $\mbf N_i$ such that for $\mbf n_i > \mbf N_i$ we have
\begin{align*}
\|U_{i,\mbf n_i} - U_i\|_{L^2(\Omega_i)} \leq C_k \left[ (\max \mbf n_i)^{\alpha+\beta} + (\max \mbf
  n_i)^{2 \alpha + \gamma} \right]^i \max_{j\leq i} \|\mathcal I_n U_j -U_j\|_{H^1(\Omega_j)},
\end{align*}
where $\mathcal I_n$ is the appropriate projection for $\Omega_j$.\\
\end{theorem}

\noindent \emph{Proof: } We prove this by induction.  Since $U_{1,\mbf
  n_1} = U_{1,n_1}$ the claim follows from Theorem \ref{convergence}
for $i =1$.  Now assume the claim is true for all $j < i$.  We use
Lemma \ref{discrete-approx-prop} to show it is true for $i$.  Using the triangle inequality we have
\begin{align*}
\|U_{i,\mbf n_i} - U_i\|_{L^2(\Omega_i)} \leq \|U_{i,\mbf n_i} - U_{i,n_i}\|_{L^2(\Omega_i)} + \|U_{i} - U_{i,n_i}\|_{L^2(\Omega_i)}.
\end{align*}
Using Theorem \ref{convergence}, we bound the second term:
\begin{align*}
\|U_{i} - U_{i,n_i}\|_{L^2(\Omega_i)} \leq C n_i^{\alpha+\beta} \|\mathcal I_n U_{i} -  U_{i}\|_{H^1(\Omega_i)}.
\end{align*}
To bound the first term we use \eqref{discrete-final-bound},
\begin{align}\label{scale-shift-b1}
\|U_{i,\mbf n_i} - U_{i,n_i}\|_{L^2(\Omega_i)} \leq C n_i^{2\beta+\gamma}E(\mbf n_{i-1}).
\end{align}
$E(\mbf n_{i-1})$ is proportional to the uniform difference of
$\tilde H_i$ and its approximation obtained through the numerical
method, $\tilde H_{i,\mbf n_{i-1}}$.  By the induction hypothesis, if
$k$ is sufficiently large, Lemma \ref{discrete-approx-prop}, tells us
that this difference tends to zero as $\mbf n_{i-1} \goto \infty$, and
the use of \eqref{discrete-final-bound} is justified.  More precisely,
the Cauchy--Schwartz inequality for each $\Omega_j$, $j < i$ and
repeated triangle inequalities results in
\begin{align}\label{scale-shift-b2}
\|\tilde H_i - \tilde H_{i,\mbf n_{i-1}}\|_{L^\infty(\Omega_i)} \leq C \sum_{j=1}^{i-1} \|U_{j} - U_{j,n_j}\|_{L^2(\Omega_j)}.
\end{align}
Combining \eqref{scale-shift-b1} and \eqref{scale-shift-b2} we complete the proof. \mqed

\begin{remark}  The requirement that $k$ is large can be made more precise using Definition \ref{admissible} with $m = \max\{l(2\alpha + \gamma),l(\alpha+\beta)\}$.  There is little restriction if $(\alpha,\beta,\gamma) = (0,0,0)$.
\end{remark}

\section{Uniform Approximation}\label{section:uniform}

In this section we describe how the above results can be
used.  First, we briefly describe how to obtain an explicit asymptotic approximation.  Second, we
use the same ideas to explain how numerics can be used to provide asymptotic approximations.  The idea we continue to exploit is that the set of invertible operators between Banach spaces is open.  Before we proceed, we define two types of uniform approximation.  Let $\{U_n^\xi\}_{\xi \geq 0}$ be a sequence, depending on the parameter $\xi$, in a Banach space such that for each $\xi$, $\|U_n^\xi-U^\xi\| \goto 0$ as $n \goto \infty$ for some $U^\xi$.

\begin{definition}  We say the sequence $\{U_n^\xi\}_{\xi \geq 0}$ is \emph{weakly uniform} if for every $\epsilon > 0$ there exists a function $N(\xi): \mathbb R^+ \goto \mathbb N$ taking finitely many values such that
\begin{align*}
\|U^\xi_{N(\xi)} - U^\xi\| < \epsilon.
\end{align*}
\end{definition}

\begin{definition} We say the sequence $\{U_n^\xi\}_{\xi \geq 0}$ is \emph{strongly uniform} (or just \emph{uniform}) if for every $\epsilon > 0$ there exists $N \in \mathbb N$ such that for $n \geq N$
\begin{align*}
\|U^\xi_{n} - U^\xi\| < \epsilon.
\end{align*}
\end{definition}

The necessity for the definition of a weakly uniform sequence is mostly a technical detail, as we do not see it {arise} in practice.  To
illustrate how it can arise we give an example.

\begin{example}
Consider the sequence
\begin{align*}
\{U_n^\xi\}_{n,\xi \geq 0}  = \left\{ \sin \xi + e^{-n^2} + e^{-(\xi-n)^2} \right\}_{n,\xi \geq 0}.
\end{align*}
For fixed $\xi$, $U_n^\xi \goto \sin \xi$.  We want, for $\epsilon
> 0$, while keeping $n$ bounded
\begin{align*}
|U_n^\xi - \sin \xi| = |e^{-n^2} + e^{-(\xi-n)^2}| < \epsilon.
\end{align*}
We choose $n > \xi$ or if $\xi$ is large enough we choose $0 < n
< \xi$.  To maintain error that is uniformly less  then $\epsilon$ we
cannot choose a fixed $n$; it must vary with respect to $\xi$.  When
relating to RHPs the switch from $n > \xi$ to $0 < n
< \xi$ is related to transitioning into the asymptotic regime.
\end{example}

\subsection{Nonlinear Steepest Descent}\label{subsection:nsd}

For simplicity, assume the RHP $[G^\xi,\Gamma^\xi]$
depends on a single parameter $\xi \geq 0$.  Further, assume $\Gamma^\xi =
\alpha(\xi) \Omega + \beta(\xi)$ where $\Omega$ is a fixed contour.  It follows that the matrix-valued functions
\begin{align*}
H^\xi(k) = G^\xi(\alpha(\xi) \Omega + \beta(\xi)),
\end{align*}
are defined on $\Omega$.  Assume there exists a matrix $H$
such that $\|H-H^\xi\|_{L^\infty(\Omega)} \goto 0$ as $\xi \goto \infty$.
Use $\Phi^\xi$ and $\Phi$ to denote the solutions of $[H^\xi;\Omega]$
and $[H;\Omega]$, respectively.  Applying  Theorem \ref{approx-prop},  $\Phi^\xi$ approximates $\Phi$ in the limit.  If everything
works out, we expect to find
\begin{align}
\Phi^\xi = \Phi + \bigo(\xi^{-\nu}), ~~\nu >0.
\end{align}
Furthermore, if we find an explicit solution to $[H;\Omega]$ we can
find an explicit asymptotic formula for the RHP
$[G^\xi;\Gamma^\xi]$.  This function $\Phi$ that solves
$[H^\xi;\Omega]$ asymptotically will be referred to as a \emph{parametrix}.

We return to Example \ref{example:disjoint} and introduce a
parameter into the problem  to demonstrate a
situation in which the RHPs on each disjoint contour decouple from
each other. 

\begin{example} \label{example:disjoint-param}
Assume $\Gamma^\xi = \Gamma^\xi_1 \cup \Gamma^\xi_2$ and
\begin{align*}
\Gamma^\xi_1 &= \alpha_1(\xi)\Omega_1 + \beta_1,\\
\Gamma^\xi_2 &= \alpha_2(\xi)\Omega_2 + \beta_2, ~~ |\beta_1 - \beta_2| > 0.
\end{align*}
Assume that each $\Omega_i$ is bounded.  We consider the $L^2$ norm of the Cauchy operator applied to a
function defined on $\Gamma^\xi_1$ and evaluated on $\Gamma^\xi_2$:
$\mathcal C_{\Gamma^\xi_1} u(z)|_{\Gamma^\xi_2}$.  Explicitly,
\begin{align*}
\mathcal C_{\Gamma^\xi_1} u(z) = \frac{1}{2 \pi i} \int_{\Gamma^\xi_1}
\frac{u(s)}{s-z} dx.
\end{align*}
This is a Hilbert--Schmidt operator, and
\begin{align}
\|\mathcal C_{\Gamma^\xi_1}\|_{\mathcal L (L^2(\Gamma^\xi_1),L^2(\Gamma^\xi_2))}^2
\leq  \int_{\Gamma^\xi_1}\int_{\Gamma^\xi_2} \frac{|dx dk|}{|x-k|^2}.
\end{align}
A simple change of variables shows that
\begin{align}\label{omega-bound}
\|\mathcal C_{\Gamma^\xi_1}\|_{\mathcal L (L^2(\Gamma^\xi_1),L^2(\Gamma^\xi_2))}^2
\leq |\alpha_1(\xi)\alpha_2(\xi)| \int_{\Omega_1} \int_{\Omega_2}
\frac{|ds dy|}{|\alpha_1(\xi)s -\alpha_2(\xi)y + \beta_1-\beta_2|^2}.
\end{align}
Since the denominator in the integral in \eqref{omega-bound} is bounded
away from zero and both $\Omega_i$ are bounded, the right-hand side tends to zero
if either $\alpha_1$ or $\alpha_2$ tend to zero. 

 This argument, with
more contours, can be used to further justify  Algorithm \ref{scale-shift} in
this limit by noting that this type of Cauchy operator
(evaluation off the contour of integration) constitutes the operator
$\mathcal T$ in \S \ref{section:setting}.  We have the representation
\begin{align*}
\mathcal C [G^\xi;\Gamma^\xi] = \begin{mat} \mathcal C [G_1^\xi;\Gamma_1^\xi] & 0 \\ 0 & \mathcal C[G_2^\xi;\Gamma_2^\xi] \end{mat} + \mathcal T^\xi,
\end{align*}
where $\|{\cal T}^\xi\|_{\mathcal L(L^2(\Gamma))} \goto 0$ as $\xi \goto
\infty$. 
\end{example}

This analysis follows similarly in some cases when $\beta_i$ depends on $\xi$.  For example, when
\begin{align} \label{betas}
\inf_{t\in S} |\beta_1(\xi) -
\beta_2(\xi)| = \delta > 0.
\end{align}
  One can extend this to the case where
\eqref{betas} is not bounded away from zero but approaches zero
slower than  $a_1(\xi)a_2(\xi)$.  For simplicity we just prove
results for $\beta_i$ being constant.  

Furthermore, the norms of the inverses are related.  When each of the
$\alpha_i(\xi)$ are sufficiently small, there exists $C > 1$ such that
\begin{align}\label{norm-equiv}
\frac{1}{C} \|\mathcal C[G;\Gamma]^{-1}\|_{\mathcal L(L^2(\Gamma))} \leq \max_i \|\mathcal C[G_i;\Gamma_i]^{-1}\|_{\mathcal L(L^2(\Gamma_i))} \leq C \|\mathcal C[G;\Gamma]^{-1}\|_{\mathcal L(L^2(\Gamma))}.
\end{align}
 Due to the simplicity of the scalings we allow the norms of the operators $\mathcal C[G_i;\Gamma_i]$ and $\mathcal C[G_i;\Gamma_i]^{-1}$ to coincide with their scaled counterparts $\mathcal C[H_i;\Omega_i]$ and $\mathcal C[H_i;\Omega_i]^{-1}$.

The choice of these scaling
parameters is not a trivial task.  We use the following rule of thumb:
\begin{assumption} \label{scalings} If the jump matrix $G$ has a factor
  $e^{\xi\theta}$ and $\beta_k$ corresponds to a $q$th order
  stationary point $\beta_k$ (\emph{i.e.}, $\theta(z) \backsim
  C(z-\beta_k)^q$), then the scaling which achieves asymptotic
  stability is $\alpha_k(\xi) = |\xi|^{-1/q}$.
\end{assumption}

We prove the validity of this assumption for the deformations
below on a case-by-case basis.  We need one final result to
guide deformations.  We start with a RHP posed on an
unbounded and connected contour.  In all cases we need to justify the truncation
of this contour, hopefully turning it into disconnected contours.

\begin{lemma} \label{contour-truncation} Assume $[G;\Gamma]$ is $k$-regular. For
  every $\epsilon > 0$ there exists a function $G_\epsilon$ defined
  on $\Gamma$ and a
  bounded contour $\Gamma_\epsilon \subset \Gamma$ such that:
\begin{itemize}
\item $G_\epsilon = I$ on $\Gamma \setminus \Gamma_\epsilon$,
\item $\|G_\epsilon- G\|_{W^{k,\infty}(\Gamma) \cap H^k(\Gamma)} < \epsilon$
\item $[G_\epsilon,\Gamma_\epsilon]$ is $k$-regular and
\begin{align*}
\|\mathcal C[G;\Gamma] - \mathcal C[G_\epsilon,\Gamma]\|_{\mathcal
  L(L^2(\Gamma))} < \epsilon\|\mathcal C_\Gamma^- \|_{\mathcal L(L^2(\Gamma))}.
\end{align*}
\end{itemize}
\end{lemma}
\emph{Proof:} A matrix-valued function $f$ is chosen  such that
\begin{itemize}
\item $f|_{\Gamma_i} \in C^\infty(\Gamma_i)$,
\item $f = I$ in a neighborhood of all intersection points,
\item $f$ has compact support, and
\item $\|(G-I)f+I - G\|_{W^{k,\infty}(\Gamma) \cap H^k(\Gamma)} < \epsilon$.
\end{itemize}
 Equate $G_\epsilon = (G-I)f + I$.  The last property follows
 immediately. \mqed

This justifies the truncation of infinite contours to finite ones and it shows this process preserves smoothness.
When it comes to numerical computations, we truncate contours when the
jump matrix is, to machine epsilon, the identity matrix.  In what
follows we assume this truncation is performed and we ignore the error induced.

\subsection{Direct Estimates}

As before, we are assuming we have a RHP that depends on a
parameter $\xi$, $[G^\xi;\Gamma^\xi]$ and $\Gamma^\xi$ is bounded.  Here we use \emph{a priori}
bounds on the solution of the associated SIE which are uniform in
$\xi$ to prove the uniform approximation.  In general, when this is possible, it is the simplest way to proceed.

Our main tool is Corollary \ref{regularity-cor}. We can easily estimate the regularity of the solution of each problem $[H_\gmi;\Omega_\gmi]$ provided we have some information about $\|\mathcal
C[G_\gmi^\xi;\Gamma_\gmi^\xi]^{-1}\|_{\mathcal L(L^2(\Gamma_\gmi))}$ or equivalently $\|\mathcal C[H^\xi_\gmi;\Omega_\gmi]^{-1}\|_{\mathcal L(L^2(\Omega_\gmi))}$.  We address how to estimate this later in this section.  First, we need a statement about how regularity is preserved throughout Algorithm \ref{scale-shift}.  Specifically, we use information from the scaled jumps $H_\gmi$ and the local inverses $\mathcal C[H_\gmi;\Omega_\gmi]^{-1}$ to estimate global regularity.  The following theorem uses this to prove strong uniformity of the numerical method.

\begin{theorem} \label{direct-uniformity} Assume
\begin{itemize}
\item  $\{[G^\xi,\Gamma^\xi]\}_{\xi \geq 0}$ is a  sequence of $k$-regular RHPs,
\item the norm of $\mathcal C[H_\gmi^\xi,\Omega_\gmi]^{-1}$ is  uniformly bounded in $\xi$,
\item $\|H_\gmi^\xi\|_{W^{k,\infty}(\Omega_\gmi)} \leq C$, and
\item $\alpha_\gmi(\xi) \goto 0$ as $\xi \goto \infty$.
\end{itemize}
Then if $k$ and $\xi$ are sufficiently large
\begin{itemize}
\item  Algorithm \ref{scale-shift} applied to $\{[G^\xi,\Gamma^\xi]\}_{\xi \geq 0}$ has solutions at each stage,
\item $\|U^\xi_j\|_{H^k(\Omega_i)} \leq P_k$ where $P_k$ depends on $\|H^\xi_i\|_{H^k(\Omega_i) \cap W^{k,\infty}(\Omega_i)}$, $\| \mathcal C[H^\xi_i;\Omega_i]^{-1}\|_{\mathcal L(L^2(\Omega_i))}$ and $\|U^\xi_j\|_{L^2(\Omega_j)}$ for $j < i$  and
\item  the approximation $U_{i,\mbf{n_i}}^\xi$ of $U_{i}^\xi$ (the solution of the SIE) at each step in Algorithm \ref{scale-shift} converges uniformly in $\xi$ as $\mbf n_i \goto \infty$, that is, convergence is strongly uniform. 
\end{itemize}
\end{theorem}

\noindent \emph{Proof: } First, we note that since $\alpha_i(\xi)
\goto 0$, \eqref{omega-bound} shows that jump matrix $\tilde H_i^\xi$
for the RHP solved at stage $i$ in Algorithm \ref{scale-shift} tends
uniformly to $H_i^\xi$.  This implies the solvability of the RHPs at each stage in Algorithm \ref{scale-shift}, and the bound
\begin{align*}
\|\mathcal C[\tilde H_i^\xi;\Omega_i]^{-1}\|_{\mathcal L(L^2(\Gamma))} \leq C \|{\cal C}[H_i^\xi;\Omega_i]^{-1}\|_{\mathcal L(L^2(\Gamma))},
\end{align*}
for sufficiently large $\xi$. As before, $C$ can be taken to be independent of $\xi$. We claim that $\|U^\xi_i\|_{H^k(\Omega_i)}$ is uniformly bounded.  We
  prove this by induction.  When $i = 1$, $U_1^\xi = \mathcal C[H_i^\xi,\Omega_i]^{-1}(H_i^\xi-I)$ and the claim follows from Corollary \ref{regularity-cor}.  Now assume the claim is true for $j< i$. All
  derivatives of the jump matrix $\tilde H^\xi_i$ depend on the Cauchy
  transform of $U_j^\xi$ evaluated away from $\Omega_j$ and $H_i^\xi$.  The former is
  bounded by the induction hypothesis and the later is bounded by
  assumption. Again, sing Corollary \ref{regularity-cor} we obtain the
  uniform boundedness of $\|U^\xi_i\|_{H^k(\Omega_i)}$. Theorem
  \ref{algorithm-convergence} implies that convergence is uniform in $\xi$.
  \mqed

The most difficult part about verifying the hypotheses of this theorem is establishing an estimate of $\|\mathcal
C[H_i^\xi;\Omega_i^\xi]^{-1}\|_{\mathcal L(L^2(\Omega_i))}$ as a function of
$\xi$.  A very useful fact is that once the solution $\Psi^\xi$ of the RHP $[G^\xi;\Gamma^\xi]$ is known then the inverse of the operator is also known:
\begin{align}\label{explicit-inverse}
\mathcal C[G^\xi;\Gamma^\xi]^{-1}u = \mathcal C_{\Gamma^\xi}^+[ u (\Psi^\xi)^+][(\Psi^\xi)^{-1}]^+ -\mathcal C_{\Gamma^\xi}^-[ u (\Psi^\xi)^+][(\Psi^\xi)^{-1}]^-.
\end{align}
This is verified by direct computation. When $\Psi^\xi$ is known
approximately, \emph{i.e.}, when a parametrix is known, then estimates
on the boundedness of the inverse can be reduced to studying the
$L^\infty$ norm of the parametrix. Then \eqref{norm-equiv} can be used to relate this to each $\mathcal C[G_i^\xi;\Gamma_i^\xi]^{-1}$ which gives an estimate for the norm of $\mathcal C[H_i^\xi;\Omega_i]^{-1}$.  We study this further
in section \ref{section:pII}.

\subsection{Failure of Direct Estimates}\label{subsection:failure}

We study a toy  RHP to motivate where the direct estimates can
fail.  Let $\phi(x)$ be a smooth function with compact support in $(-1,1)$ satisfying $\max_{[-1,1]} |\phi(x)| =1/2 $.  Consider the following scalar RHP for a function $\mu$:
\begin{align}\label{mu-growth}
\mu^+(x) &= \mu^-(x) (1 + \phi(x)(1 + \xi^{-1/2} e^{i\xi x})), \\
\mu(\infty) &= 1, ~~x \in [-1,1],~~ \xi > 0.
\end{align}

This problem can be solved explicitly but we study it from the linear
operator perspective instead. From the boundedness assumption on
$\phi$, a Neumann series argument gives the invertiblity  of the
singular integral operator and uniform boundedness of the $L^2$
inverse in $\xi$.  Using the estimates in Corollary
\ref{regularity-cor} we obtain useless bounds, which all grow with $\xi$.  Intuitively, the solution to \eqref{mu-growth} is close, in
$L^2$ to the solution to
\begin{align}\label{mu-limit}
\nu^+(x) &= \nu^-(x) (1 + \phi(x)),\\
\nu(\infty) &= 1,~~ x \in [-1,1],~~ \xi > 0,
\end{align}
which trivially has uniform bounds on its Sobolev norms.  In the next
section we introduce the idea of a numerical parametrix which helps resolve this complication.

\subsection{Extension to Indirect Estimates}


In this section we assume minimal hypotheses for dependence of the
sequence $\{[G^\xi;\Gamma^\xi]\}_{\xi \geq 0}$ on $\xi$.  Specifically
we require only  that the map $\xi \mapsto H_i^\xi$ is continuous from
$\mathbb R^+$ to $L^\infty(\Omega_i)$ for each $i$.  We do not want to
hypothesize more as that would alter the connection to the method of
nonlinear steepest descent which only requires uniform convergence of
the jump matrix.  In specific cases, stronger results can be obtained
by requiring the map $\xi \mapsto H_i^\xi$ to be continuous from
$\mathbb R^+$ to $W^{k,\infty}(\Omega_i)$.  

The fundamental result we need to prove a uniform approximation
theorem is the continuity of Algorithm \ref{scale-shift} with respect to uniform perturbations in the jump matrix.  With the jump matrix $G$ we associated $H_j$, the scaled restriction of $G$ to $\Gamma_j$.  With $G$ we also associated $U_j$, the solution of the SIE obtained from $[\tilde H_j;\Omega_j]$.  In what follows we will have another jump matrix $J$ and analogously we use $K_j$ to denote the scaled restriction of $J$ and $P_j$ to denote the solution of the SIE obtained from $[\tilde K_j,\Omega_j]$.

\begin{lemma} \label{algorithm-continuity} Assume
  $\{[G^\xi,\Gamma^\xi]\}_{\xi \geq 0}$ is a sequence of $1$-regular
  RHPs such that $\xi \mapsto H_i^\xi$ is continuous from
  $\mathbb R^+$ to $L^\infty(\Omega_i)$ for each $i$.  Then for
  sufficiently large, but fixed,
  $\mbf n_i$ the map $\xi \mapsto U_{i,\mbf n_i}^\xi$ is continuous from $\mathbb R^+$ to $L^2(\Omega_i)$ for each $i$.
\end{lemma}

\noindent \emph{Proof: }  We prove this by induction on $i$.  For $i =
1$ the claim follows from Lemma \ref{discrete-approx-prop}.  Now
assume the claim is true for $j<i$ and we prove it holds for $i$.  We
show the map is continuous at $\eta$ for $\eta \geq 0$.  First, from
Lemma \ref{discrete-approx-prop}
\begin{align*}
\|U_{i,\mbf n_i}^\eta - U_{i,\mbf n_i}^\xi\|_{L^2(\Omega)} \leq C n_i^{2\alpha+\gamma} E(\xi,\eta),
\end{align*}
where $E(\xi,\eta)$ is proportional to $\|\tilde H_{i,\mbf n_{i-1}}^\eta -\tilde H_{i,\mbf n_{i-1}}^\xi\|_{L^\infty(\Omega_i)}$.  A similar argument as in Theorem \ref{algorithm-convergence} gives
\begin{align*}
\|\tilde H_{i,\mbf n_{i-1}}^\eta -\tilde H_{i,\mbf
  n_{i-1}}^\xi\|_{L^\infty(\Omega_i)} \leq C(\eta,\mbf n_i) \sum_{j=1}^{i-1} \|U^\eta_{j,\mbf n_j} - U^\xi_{j,\mbf n_j}\|_{L^2(\Omega_j)},
\end{align*}
for $|\xi-\eta|$ sufficiently small.  By assumption, the right-hand
side tends to zero as $\xi \goto \eta$ and this proves the lemma. \mqed

It is worthwhile noting that the arguments in Lemma
\ref{algorithm-continuity} show the same continuity for the
infinite-dimensional, non-discretized problem.  Now we show weak uniform convergence of the numerical scheme on compact sets.

\begin{proposition} \label{compact-uniformity} Assume
  $\{[G^\xi,\Gamma^\xi]\}_{\xi \geq 0}$ is a sequence of $k$-regular
  RHPs such that all the operators in Algorithm
  \ref{scale-shift} are invertible for every $\xi$.  Assume that $k$
  is sufficiently large so that the approximations from Algorithm \ref{scale-shift} converge for every $\xi \geq 0$. Then there exists a vector-valued function  $\mbf N(i,\xi)$ that takes finitely many values such that
\begin{align*}
\|U^\xi_{i,\mbf N(i,\xi)} - U_i^\xi\|_{L^2(\Omega_i)} < \epsilon.
\end{align*}
Moreover if the numerical method is of type $(0,0,0)$ then convergence is strongly uniform.
\end{proposition}

\noindent \emph{Proof: } Let $S \subset \mathbb R^+$ be compact.  It follows from Lemma \ref{algorithm-continuity} that the  function $E(\xi,\mbf n,i) = \|U^\xi_{i,\mbf n} -
U_{i}^\xi\|_{L^2(\Gamma_i)}$ is a continuous function of $\xi$ for fixed $\mbf n$. For $\epsilon > 0$ find an
$\mbf n_\xi$ such that  $E(\xi,\mbf n_\xi,i) < \epsilon/2$.  By continuity, define
$\delta_\xi(\mbf n_\xi) > 0$ so that $E(s,\mbf n_\xi,i) < \epsilon$ for $|s-\xi| <
\delta_\xi$. Define the ball $B(\xi,\delta) = \{s \in \mathbb R^+: |s-\xi| < \delta\}$.  The open sets $\{B(\xi,\delta_\xi)\}_{ \xi \in  S}$ cover $S$ and we can select
a finite subcover $\{B(\xi_j,\delta_{\xi_j})\}_{j=1}^N$.  We have
$E(s,\mbf n_{\xi_j},i) < \epsilon$ whenever $s \in
B(\xi_j,\delta_{\xi_j})$. To prove the claim for a method of type
$(0,0,0)$, we use the fact that $\delta_\xi$ can be taken independent of
$\mbf n_\xi$ and that $E(s,\mbf n,i) < \epsilon$ for every $\mbf n >
\mbf n_\xi$. \mqed

\begin{definition} \label{admissible-para} Given a sequence of  $k$-regular RHPs $\{[G^\xi,\Gamma^\xi]\}_{\xi \geq 0}$ such that
\begin{itemize}
\item $\Gamma^\xi = \Gamma_{1}^\xi \cup \cdots \cup \Gamma_{\ell}^\xi$, and 
\item $\Gamma_\gmi^\xi = \alpha_\gmi(\xi) \Omega_\gmi + \beta_\gmi$,
\end{itemize}
another sequence of $k$-regular RHPs $\{[J^\xi,\Sigma^\xi]\}_{\xi\geq0}$ is said to be a numerical parametrix if
\begin{itemize}
\item $\Sigma^\xi = \Sigma_{1}^\xi \cup \cdots \cup \Sigma_\ell^\xi$,
\item $\Sigma_\gmi^\xi = \gamma_\gmi(\xi) \Omega_\gmi + \sigma_\gmi$,
\item For all $i$
\begin{align}
J^\xi(\gamma_i(\xi) k + \sigma_i) - G^\xi(\alpha_i(\xi) k + \beta_i) \goto 0,
\end{align}
uniformly on $\Omega_i$ as $\xi \goto \infty$,
\item the norms of the  operators and inverse operators at each step in Algorithm \ref{scale-shift} are uniformly bounded in $\xi$, implying uniform boundedness of $J^\xi$ in $\xi$, and
\item the approximation $P_{i,\mbf{n_i}}^\xi$ of $P_{i}^\xi$ (the solution of the SIE) at each step in Algorithm \ref{scale-shift} converges uniformly as $\min \mbf n_i \goto \infty$. 
\end{itemize}
\end{definition}

This definition hypothesizes {desirable conditions on a nearby limit problem} for the
sequence $\{[G^\xi,\Gamma^\xi]\}_{\xi \geq 0}$. Under the assumption of this nearby limit problem we are able to obtain a uniform approximation for the solution of the original RHP.

\begin{lemma}  \label{parametrix-uniformity}  Assume there exists a
  numerical parametrix $\{J_\xi,\Sigma_\xi\}_{\xi>0}$  for a sequence of RHPs $\{[G_\xi,\Gamma_\xi]\}_{\xi \geq 0}$.  Then for every $\epsilon > 0$ there exists  $\mbf N_i$ and $T>0$  such that, at each stage in Algorithm \ref{scale-shift},
\begin{align}\label{uniform-state}
\|U^\xi_{i,\mbf N_i} - U^\xi_{i}\|_{L^2(\Omega_i)} < \epsilon \mbox{ for } \xi > T.
\end{align}
Furthermore, if the numerical method is of type $(0,0,0)$, then
\eqref{uniform-state} is true with $\mbf N_i$ replaced by any $\mbf M_i > \mbf N_i$.

\end{lemma}

\noindent \emph{Proof: }  At each stage in Algorithm \ref{scale-shift} we have
\begin{align}\label{SIE-sol-bound}
\|U^\xi_{i,\mbf n_i} - U^\xi_i\|_{L^2(\Omega_i)} \leq \|U^\xi_{i,\mbf n_i} - P^\xi_{i,\mbf n_i}\|_{L^2(\Omega_i)} + \|P^\xi_{i,\mbf n_i} - P_i^\xi\|_{L^2(\Omega_i)} + \|P_i^\xi - U_i^\xi\|_{L^2(\Omega)}.
\end{align}
Since $P_{i,\mbf n_i}^\xi$ originates from a numerical parametrix we
know that $\|P^\xi_{i,\mbf n_i} - P_i^\xi\|_{L^2(\Omega_i)} \goto 0$
uniformly in $\xi$ as $\mbf n_i$ is increased.  Furthermore,
$\|P_i^\xi - U_i^\xi\|_{L^2(\Omega)}$ depends only on $\xi$ and tends
to zero as $\xi \goto \infty$.  The main complication comes from the fact
that a bound on $\|U^\xi_{i,\mbf n_i} - P^\xi_{i,\mbf
  n_i}\|_{L^2(\Omega_i)}$ from \eqref{discrete-bound}  depends on both
$\mbf n_{i-1}$ and $\xi$ if the method is not of type $(0,0,0)$.  The same arguments as in Lemma
\ref{algorithm-continuity} show this tends to zero.  Therefore we choose $\mbf n_i$ large
enough so that the second term in \eqref{SIE-sol-bound} is less than
$\epsilon/3$.  Next, we choose $\xi$ large enough so that the sum of
the remaining terms is less than $2/3 \epsilon$.  If the method is of
type $(0,0,0)$ this sum remains less than $\epsilon$ if $\mbf n_i$ is
replaced with $\mbf n$ for $\mbf n > \mbf n_i$.  This proves the
claims. \mqed

Now we prove the uniform approximation theorem.

\begin{theorem} \label{uniform-thm} Assume
  $\{[G^\xi,\Gamma^\xi]\}_{\xi \geq 0}$ is a sequence of $k$-regular
  RHPs for $k$ sufficiently large so that Algorithm
  \ref{scale-shift} converges for each $\xi$.  Assume there exists a
  numerical parametrix as $\xi \goto \infty$.  Then Algorithm
  \ref{scale-shift} produces a weakly uniform approximation to the
  solution of $\{[G^\xi,\Gamma^\xi]\}_{\xi \geq 0}$.  Moreover, convergence is strongly uniform if the method is of type $(0,0,0)$.
  \end{theorem}

\noindent \emph{Proof: }  Lemma \ref{parametrix-uniformity} provides an $M > 0$ and $\mbf N_1(i)$ such that if $\xi > M$ then
\begin{align*}
\|U_{i,\mbf N_1(i)}^\xi - U_i^\xi\|_{L^2(\Omega_i)} < \epsilon, \mbox{ for every } i.
\end{align*}
According to  Theorem \ref{compact-uniformity} there is an $\mbf N_2(\xi,i)$ such that
\begin{align*}
\|U_{i,\mbf N_2(\xi,i)}^\xi - U_i^\xi\|_{L^2(\Omega_i)} < \epsilon,  \mbox{ for every } i.
\end{align*}
The function
\begin{align*}
\mbf N(\xi,i) = \begin{choices} \mbf N_1(i), \when \xi > M,\\
N_2(\xi,i), \when \xi \leq \mbf M, \end{choices}
\end{align*}
satisfies the required properties for weak uniformity.  Strong
uniformity follows in a similar way from Lemma
\ref{parametrix-uniformity} and Theorem \ref{compact-uniformity}. \mqed

\begin{remark}  This proves weak uniform convergence of the numerical method for the toy problem introduced in \S \ref{subsection:failure}:  we can take the RHP for $\nu$ as a numerical parametrix.
\end{remark}

The seemingly odd restrictions for the general theorem are a
consequence of poorer operator convergence rates when $n$ is large.  A
well-conditioned numerical method does not suffer from this issue.  It
is worth noting that using direct estimates is equivalent to {requiring} that the RHP itself satisfies the properties of a numerical parametrix.

In what follows, we want to show a given sequence of RHPs is a
numerical parametrix.  The reasoning for the following result is
two-fold. First, we hypothesize only conditions which are easily
checked in practice.  Second, we want to connect the stability of numerical approximation with the use of local, model problems in nonlinear steepest descent.

\begin{proposition} \label{local-prop} Assume
\begin{itemize}
\item  $\{[J^\xi,\Sigma^\xi]\}_{\xi \geq 0}$ is a sequence of $k$-regular RHPs,
\item $\mathcal C[K_i^\xi,\Omega_i]^{-1}$ has norm uniformly bounded in $\xi$,
\item $\|K_i^\xi\|_{W^{k,\infty}(\Omega_i)} \leq C$, and
\item $\gamma_i(\xi) \goto 0$ as $\xi \goto \infty$.
\end{itemize}
Then if $k$ and $\xi$ are sufficiently large
\begin{itemize}
\item  Algorithm \ref{scale-shift} applied to $\{[J^\xi,\Sigma^\xi]\}_{\xi \geq 0}$ has solutions at each stage and
\item  $\{[J^\xi,\Sigma^\xi]\}_{\xi \geq 0}$ satisfies the last two
  properties of a numerical parametrix (Definition \ref{admissible-para}).
\end{itemize}
\end{proposition}

\noindent \emph{Proof: } The proof of this is essentially the same as
Theorem \ref{direct-uniformity} \mqed.

\begin{remark}
Due to the decay of $\gamma_i$, the
invertiblity of each of $\mathcal C [K_i^\xi;\Omega_i]$ is equivalent to that
of $\mathcal C[G^\xi;\Gamma^\xi]$.
\end{remark}

This proposition states that a numerical parametrix only needs to be locally reliable; we can consider each shrinking contour as a separate RHP as far as the analysis is concerned.

\section{A Numerical Realization} \label{section:realize}


In \cite{SORHFramework}, a numerical framework was constructed for
computing solutions to RHPs, based on a method used in
\cite{SOPainleveII} for computing solutions to the undeformed
Painlev\'e II RHP.  This framework is based on Chebyshev
interpolants.  Consider the RHP $[G;\Gamma]$, $\Gamma =
\Gamma_1 \cup \cdots \cup \Gamma_\ell$,  where each $\Gamma_\gmi$ is bounded
and is a
M\"obius transformation of the unit interval:
\begin{align*}
M_\gmi([-1,1]) = \Gamma_\gmi.
\end{align*}

\begin{definition} \label{Chebypoints} The Chebyshev points of the second kind are
\begin{align*}
\mbf x^{[-1,1],n} = \begin{mat} x_1^{[-1,1],n} \\
  \vdots \\ x_n^{[-1,1],n} \end{mat} = \begin{mat} -1 \\ \cos \pi
  \left(1-\frac{1}{n-1}\right) \\ \vdots \\ \cos \frac{\pi}{n-1} \\
  1 \end{mat}.
\end{align*}
\end{definition}

The mapped Chebyshev points are denoted
\begin{align*}
\mbf x^{\gmi,n} = M_\gmi(\mbf x^{[-1,1]}).
\end{align*}

Given a continuous function $f_\gmi$ defined on $\Gamma_\gmi$ we can find a unique interpolant at $\mbf x^{\gmi,n}$ using mapped Chebyshev polynomials.
Given a function, $f$ defined on the whole of $\Gamma$, we define $\mathcal I_n$ to be this interpolation projection applied
to the restriction of $f$ on each $\Gamma_\gmi$.  Clearly,
\begin{align*}
\mathcal I_n : H^1(\Gamma) \goto H^1(\Gamma)
\end{align*}
{and, because $\mbf x^{j,n}$ contains all junction points,}
\begin{align*}
\mathcal I_n : H^1_z(\Gamma) \goto H^1_z(\Gamma).
\end{align*}

The framework in \cite{SORHFramework} is given by the pair $(\mathcal I_n, \mathcal
I_n)$ and the matrix $\mathcal C_n[G;\Gamma]$ is equal to $\mathcal
I_n \mathcal C[G;\Gamma]\mathcal I_n$ with some unbounded components
subtracted; obeying the requirement that the two operators agree on $H^1_z(\Gamma)$.

Now we address the properties required in Definition
\ref{admissible}.

\begin{lemma} When $G \in W^{1,\infty}(\Gamma)$, the numerical method in \cite{SORHFramework} satisfies:
\begin{itemize}
\item  $\mathcal I_n$ is uniformly bounded in $n$ from $C(\Gamma)$ to $L^2(\Gamma)$ when $\Gamma$ is bounded.
\item $\|\mathcal C_n[G;\Gamma]\|_{\mathcal L(H^1_z(\Gamma),Y_n)} \leq C(1 + \Norms)$.
\item $\|\mathcal T_n[G;\Gamma]\|_{\mathcal L(X_n,Y_n)} \leq C n^2 \Norms$.
\item $\|\mathcal I_nu-u\|_{H^{1}(\Gamma)} \leq C_s n^{2-s} \|u\|_{H^s(\Gamma)}$.
\end{itemize}

\end{lemma}
\noindent \emph{Proof: }  First, note that these constants depend on $\Gamma$. Using the Dirichlet kernel one proves that $\mathcal I_n$ is uniformly bounded from $C(\Gamma)$ to an $L^2$ space with the Chebyshev weight \cite{atkinson}.  The norm on this weighted space dominates the usual $L^2(\Gamma)$ norm, proving the first result.  For the second statement we take $u \in H^1_z(\Gamma)$ and consider
\begin{align*}
\|\mathcal I_n - \mathcal I_n (G-I) \mathcal C_\Gamma^- u\|_{L^2(\Gamma)} \leq \|\mathcal I_n\|_{\mathcal L(C(\Gamma),L^2(\Gamma))} ( 1 +\|G-I\|_{L^\infty(\Gamma)} \|\mathcal C_\Gamma^-\|_{\mathcal L(H^1_z(\Gamma),H^1(\Gamma))} \|u\|_{H^1(\Gamma)}).
\end{align*}
Since $Y_n$ is equipped with the $L^2(\Gamma)$ norm and $\|\mathcal C_\Gamma^-\|_{\mathcal L(H^1_z(\Gamma),H^1(\Gamma))} = \|\mathcal C_\Gamma^-\|_{\mathcal L(L^2(\Gamma))}$ we obtain the second property.  We then use for $u \in X_n$ (see \S \ref{section:setting}) that $\|u\|_{H^1(\Gamma)} \leq Cn^2 \|u\|_{L^2(\Gamma)}$ to obtain
\begin{align*}
\|\mathcal T_n[G;\Gamma]\|_{\mathcal L(X_n,Y_n)} \leq C n^2 \Norms.
\end{align*}
The last statement follows from estimates in \cite{quarteroni} for
the pseudo-spectral derivative. \mqed

 The final property we need to obtain an admissible numerical method, the
boundedness of the inverse, is a very difficult problem.  We can
verify, \emph{a posteriori}, that the norm of the inverse does not grow
too much.  In general, for this method, we see at most logarithmic
growth.  We make the following assumption.

\begin{assumption} \label{bounded-assumption}For the
  framework in \cite{SORHFramework} we assume that whenever $[G;\Gamma]$ is $1$-regular and $\mathcal
  C[G;\Gamma]^{-1}$ exists on $L^2(\Gamma)$ as a bounded operator we have
  for $n > N$
\begin{align}\label{op-bound}
\|\mathcal C_n [G;\Gamma]^{-1} \|_{\mathcal L(Y_n,X_n)}  \leq C n^\beta \|\mathcal C[G;\Gamma]^{-1}\|_{\mathcal L(L^2(\Gamma))}, ~~\beta >
0.
\end{align}
\end{assumption}

\begin{remark} This assumption is known to be true for a similar collocation method on the unit circle using Laurent monomials \cite{prossdorf}.
\end{remark}

With this assumption the numerical method associated with $(\mathcal I_n,\mathcal I_n)$ is of type $(0,\beta,2)$.  In light of Theorem \ref{convergence}, we expect spectral convergence
and the bound in Assumption \ref{bounded-assumption} does not prevent
convergence.  We combine Assumption \ref{bounded-assumption}, Theorem \ref{convergence} and Theorem \ref{regularity}  to obtain
\begin{align}\label{framework-estimate}
\|u - u_n\|_{L^2(\Gamma)} \leq  C(\|\mathcal C[G;\Gamma]^{-1}\|_{\mathcal L(L^2(\Gamma))} (1 + \Norms ) n^{2+\beta-k} \|u\|_{H^k(\Gamma)}.
\end{align}

\section{Application to Painlev\'e II}\label{section:pII}

\begin{figure}[ht]
\centering
\includegraphics[width=.8\linewidth]{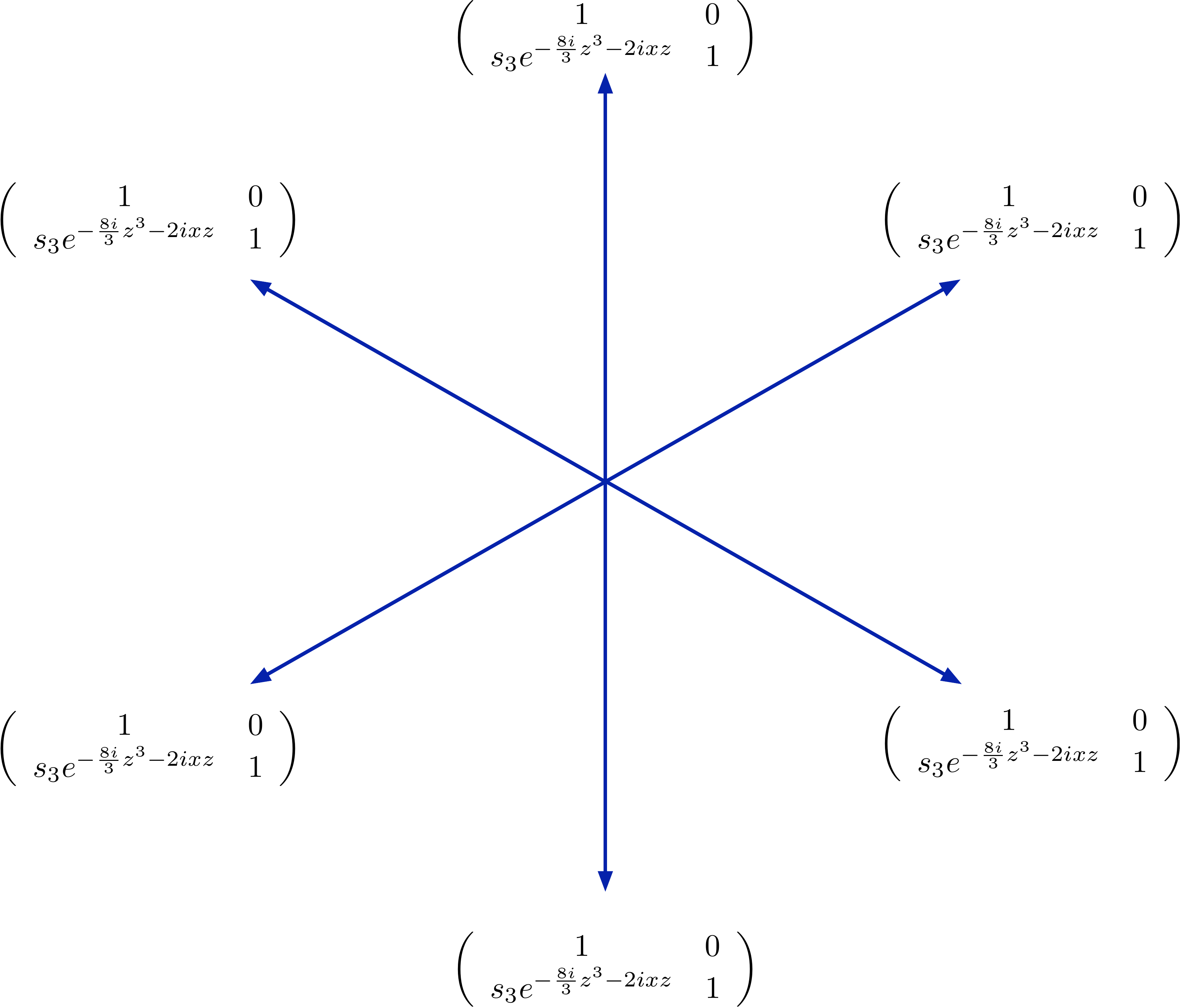}
\caption{\label{UndeformedContour} The contour and jump matrix for the \PII RHP.}
\end{figure}

We present the RHP for the solution of the Painlev\'e II ODE \eqref{PII}.
Let $\Gamma = \Gamma_1 \cup \cdots \cup \Gamma_6$ with $\Gamma_\gmi = \{ s
e^{i\pi (\gmi/3-1/6)}: s \in \mathbb R^+\}$, \emph{i.e.}, $\Gamma$
consists of six rays emanating from the origin, see Figure 
\ref{UndeformedContour}.  The jump matrix is defined by
$G(\lambda) = G_\gmi(\lambda)$ for $\lambda \in \Gamma_\gmi$, where 
\begin{align*}
G_\gmi(x;\lambda) = G_\gmi(\lambda) = \begin{choices} \begin{mat} 1
    & s_\gmi e^{-i8/3 \lambda^3 -2i x \lambda} \\0 & 1 \end{mat} \when
  \gmi \mbox{ is even},\\
\begin{mat} 1 & 0 \\  s_\gmi e^{i8/3 \lambda^3 +
    2ix\lambda} & 1 \end{mat} \when \gmi \mbox{ is odd}. \end{choices}
\end{align*}
From the solution $\Phi$ of $[G;\Gamma]$, the Painlev\'e function is recovered by the formula
\begin{align*}
u(x) = \lim_{z \goto \infty} z \Phi_{12}(z),
\end{align*}
where the subscripts denote the $(1,2)$ entry.  This RHP was solved numerically in
\cite{SOPainleveII}.

\begin{figure}[ht]
\centering
\includegraphics[width=.9\linewidth]{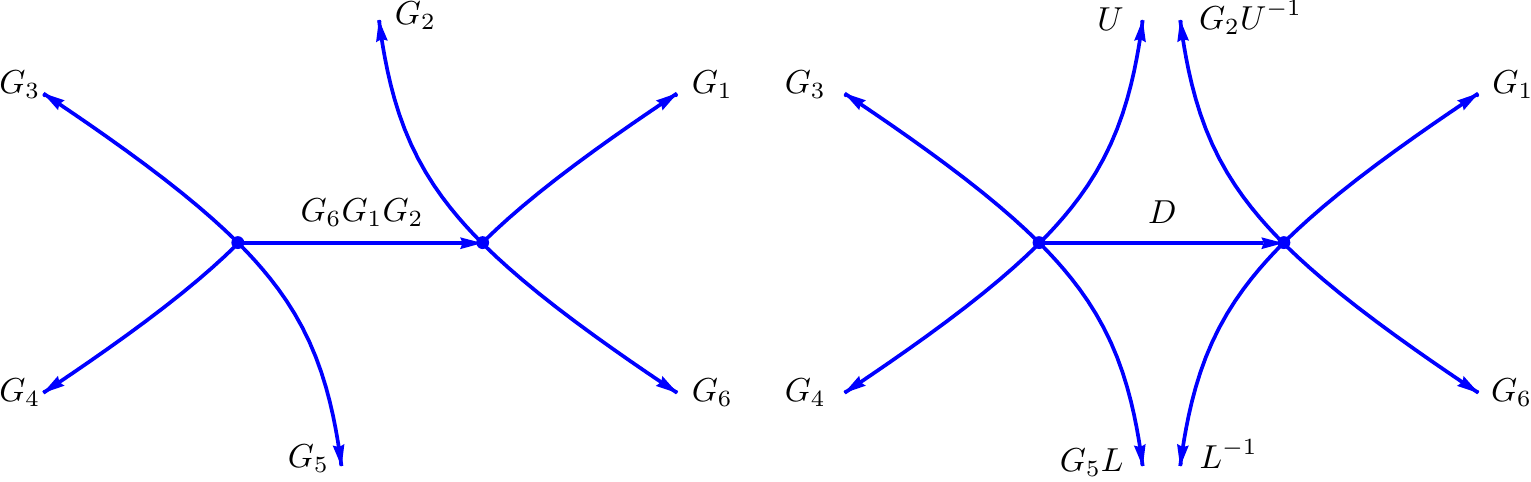}
\caption{Left: Initial deformation along the paths of steepest descent. Right: The deformed contour after lensing. \label{InitialDeformedContour}}
\end{figure}

For large $|x|$, the jump matrices $G$ are increasingly
oscillatory.  We combat this issue by deforming the contour so
that these oscillations turn to exponential decay.  To simplify this
procedure, and to start to mold the RHP into the abstract form in
\S \ref{section:uniform}, we first rescale the RHP.  If we let
$z = \sqrt{|x|} \lambda$, then the jump contour $\Gamma$ is
unchanged, and
\begin{align*}
\Phi^+(z) = \Phi^+(\sqrt{|x|} \lambda ) = \Phi^-(\sqrt{|x|}
\lambda)G(\sqrt{|x|} \lambda) = \Phi^-(z) G(z),
\end{align*}
where $G(z) = G_\gmi(z)$ on $\Gamma_\gmi$ with
\begin{align*}
G_\gmi(z) = \begin{choices} \begin{mat} 1 & s_\gmi e^{-\xi
      \theta(z)} \\0  & 1 \end{mat} \when \gmi \mbox{ is even},\\
\begin{mat} 1 & 0 \\ s_\gmi e^{\xi\theta(z)} & 1 \end{mat} \when
\gmi \mbox{ is odd}, \end{choices}
\end{align*}
$\xi = |x|^{3/2}$ and
\begin{align*}
\theta(z) = \frac{2i}{3} \left( 4z^3 + 2 e^{i \arg x} z\right).
\end{align*}
Then
\begin{align}\label{PII-reconstruct}
u(x) = \lim_{\lambda \goto \infty} \lambda \Phi_{12}(x;\lambda) =
\sqrt{x} \lim_{\lambda \goto \infty} z \Phi_{12}(x;z).
\end{align}
We assume that $s_1s_3 > 1$ and $x < 0$. We deform $\Gamma$ to pass through the stationary points $\pm
1/2$, resulting in the RHP on the left of Figure
\ref{InitialDeformedContour}.

\begin{figure}[ht]
\centering
\includegraphics[width=.5\linewidth]{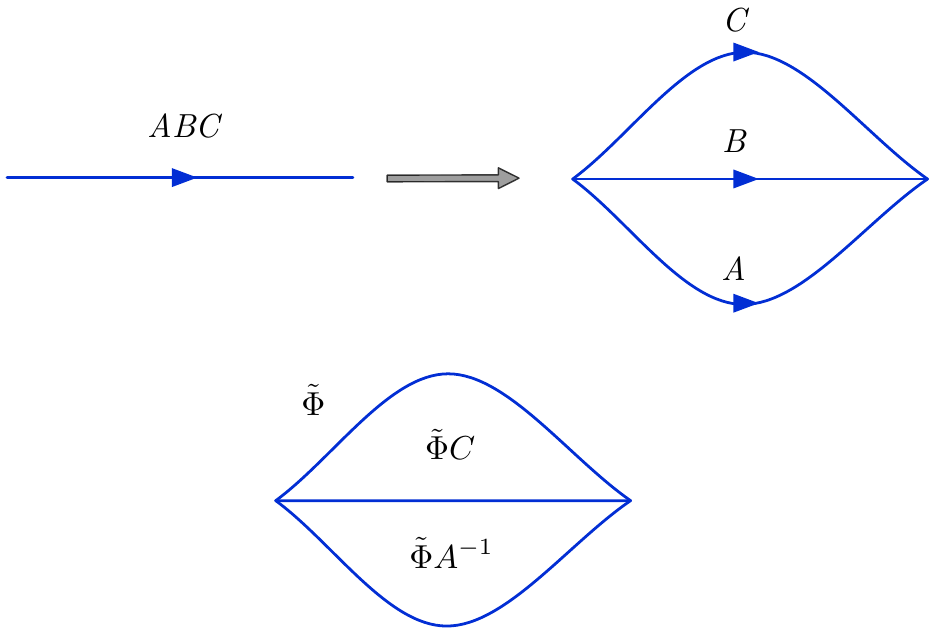}
\caption{Demonstration of the lensing process. \label{Lensing}}
\end{figure}

The function
\begin{align*}
G_6G_1G_2 = \begin{mat} 1 - s_1s_3 & s_1 e^{-\xi \theta} \\ s_2
  e^{\xi \theta} & 1 + s_1s_2 \end{mat},
\end{align*}
 has terms with $\exp(\pm \xi \theta)$. It cannot decay to the identity
when deformed in the complex plane.  We can resolve this
issue by using \emph{lensing} \cite{DeiftZhouAMS}.  Suppose we factor the jump function as $ABC$.  It is possible to separate
the single contour into three contours, as in Figure \ref{Lensing},
assuming $A$ and $C$ are analytic between the original contour and
continuous up to the new contour. If $\tilde \Phi$ satisfies the jumps
on the split contour, it is clear that we can recover $\Phi$ by
defining $\Phi = \tilde \Phi C$ between the top contour and the
original contour, $\Phi = \tilde \Phi A^{-1}$ between the original
contour and the bottom contour, and $\Phi = \tilde \Phi$ everywhere
else, see the bottom of Figure \ref{Lensing}.  As in the standard
deformation of contours, the limit at infinity is unchanged.

Now consider the $LDU$ factorization:
\begin{align*}
G_6G_1G_2 = LDU = \begin{mat} 1 & 0 \\ e^{-\zeta \theta}
  \frac{s_1}{1-s_1s_3} & 1 \end{mat} \begin{mat} 1-s_1s_3 & 0\\ 0 &
  \frac{1}{1-s_1s_3} \end{mat} \begin{mat} 1 & e^{\zeta\theta}
  \frac{s_1}{1-s_1s_3} \\  0 & 1 \end{mat}.
\end{align*}
$U$ decays near $i \infty$, $L$ decays near $-i \infty$, both go to to
the identity matrix at infinity and $D$ is constant.  Moreover, the
oscillators on $L$ and $U$ are precisely those of the original $G$
matrices.  Therefore, we reuse the path of steepest descent, and
obtain the deformation on the right of  Figure
\ref{InitialDeformedContour}.
The $LDU$ decomposition is valid under the assumption $s_1s_2 > 1$.

\subsection{Removing the  connected contour}

Although the jump matrix $D$ is non-oscillatory (in fact, constant), it
is still incompatible with the theory presented in
\S \ref{section:uniform}: we need the jump matrix to approach
the identity matrix away from the stationary points.  Therefore, it is
necessary to remove this connecting contour.  Since $D =
\diag(d_1,d_2)$ is diagonal, we can solve $P^+ = P^- D$ with
$P(\infty) = I$ on $(-1/2,1/2)$ in closed form \cite{FokasPainleve}:
\begin{align*}
P(z) = \begin{mat} \left(\frac{2x+1}{2x-1}\right)^{{i \log d_1}/{2
      \pi}} & 0 \\ 0  & \left(\frac{2x+1}{2x-1}\right)^{{i \log d_2}/{2
      \pi}} \end{mat}.
\end{align*}

\begin{figure}[ht]
\centering
\includegraphics[width=.9\linewidth]{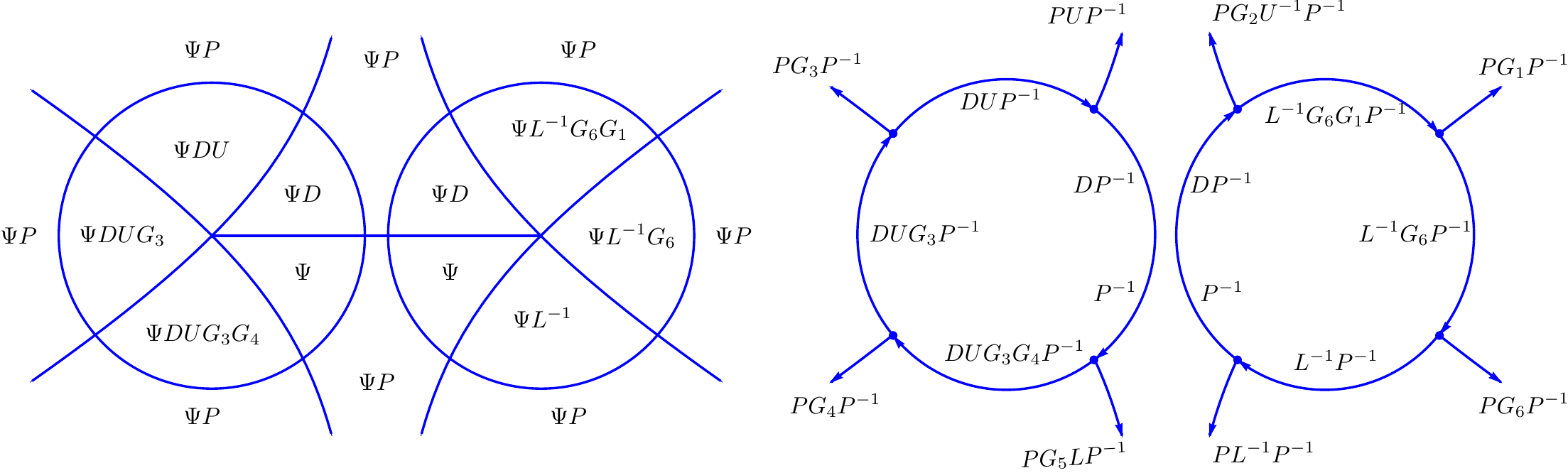}
\caption{Left: Definition of $\Phi$ in terms of $\Psi$.  Right: Jump
  contour for $\Psi$. \label{RemoveParametrixCurve}}
\end{figure}

This parametrix solves the desired RHP for any choice of
branch of the logarithm.  However, we must choose the branch so that
the singularity is square integrable \cite{FokasPainleve}. In
this case this is accomplished by choosing the standard choice of branch.

We write
\begin{align*}
\Phi = \Psi P.
\end{align*}
Since $P$ satisfies the required jump on $(-1/2,1/2)$, $\Psi$ has no
jump there.  Moreover, on each of the remaining curves we
have
\begin{align*}
\Psi^+ = \Phi^+ P^{-1} = \Phi^- G P^{-1} = \Psi^{-1} P GP^{-1},
\end{align*}
and our jump matrix becomes $PGP^{-1}$.  Unfortunately, we have
introduced singularities at $\pm 1/2$ and the theory of \S
\ref{section:uniform} requires some smoothness of the jump
matrix.   This motivates alternate definitions for $\Psi$ in circles
around the stationary points.  In particular, we define $\Phi$ in
terms of $\Psi$ by the left panel of Figure
\ref{RemoveParametrixCurve}, where $\Psi$ has the jump matrix defined
in the right.  A quick check demonstrates that this definition of
$\Phi$ indeed satisfies the required jump relations.

{We are ready to apply Algorithm~\ref{scale-shift}.  Define
	$$\Omega = \{z : \|z\| = 1\} \cup \{r e^{i \pi/4} : r \in (1,2)\} \cup \{r e^{3i \pi/4} : r \in (1,2)\} \cup \{r e^{-3 i \pi/4} : r \in (1,2)\}\cup \{r e^{-i \pi/4} : r \in (1,2)\}.$$
In accordance with Assumption \ref{scalings}, we have
	$$\Gamma_1 = {1 \over 2} + \xi^{-1/2} \Omega\quad\hbox{and}\quad\Gamma_2 = -{1 \over 2} + \xi^{-1/2} \Omega,$$
with the jump matrices defined according to Figure~\ref{RemoveParametrixCurve}.  Paths of steepest descent are now local paths of steepest descent. }


\subsection{Uniform Approximation}

We have isolated the RHP near the stationary points{, and constructed a numerical algorithm to solve the deformed RHP.  We show that this numerical algorithm approximates the true solution to the RHP.  In order to analyze the error, we  introduce the local model problem for this RHP following \cite{FokasPainleve}}.  

Define  the Wronskian matrix of parabolic cylinder functions $D_\nu(\zeta)$,
\begin{align}
Z_0(\zeta) = \begin{mat} D_{-\nu-1}(i\zeta) & D_\nu(\zeta)\\\frac{d}{d\zeta} D_{-\nu-1}(i\zeta) & \frac{d}{d\zeta} D_{\nu}(\zeta) \end{mat},
\end{align}
and the constant matrices
\begin{align*}
H_{k+2} &= e^{i\pi(\nu+1/2)\sigma_3} H_k e^{i\pi(\nu+1/2)\sigma_3}, H_0 = \begin{mat} 1 & 0 \\ h_0 & 1\end{mat}, H_2 = \begin{mat} 1 & h_1 \\ 0 & 1 \end{mat}, \sigma_3 = \begin{mat} 1 & 0 \\ 0 & -1 \end{mat}, \\
h_0 &= -i \frac{\sqrt{2\pi}}{\Gamma(\nu+1)}, ~~ h_1 = \frac{\sqrt{2 \pi}}{\Gamma(-\nu)}e^{i\pi\nu}.
\end{align*}
The sectionally holomorphic function $Z(\zeta)$ is defined as
\begin{align*}
Z(\zeta) = \begin{choices} Z_0(\zeta), \when \arg \zeta \in (-\pi/4,0),\\
Z_1(\zeta), \when \arg \zeta \in (0, \pi/2),\\
Z_2(\zeta), \when \arg \zeta \in (\pi/2,\pi),\\
Z_3(\zeta), \when \arg \zeta \in (\pi,3\pi/2),\\
Z_4(\zeta), \when \arg \zeta \in (3\pi/2,7\pi/4).
\end{choices}
\end{align*}
This is used to construct the local solutions
\begin{align*}
\hat \Psi^r(z) &= B(z)(-h_1/s_3)^{-\sigma_3/2}e^{it\sigma_3/2} 2^{-\sigma_3/2}\begin{mat}\zeta(z) & 1 \\ 1 & 0 \end{mat} Z(\zeta(z))(-h_1/s_3)^{\sigma_3/2}.\\
\hat \Psi^l(z) &= \sigma_2 \hat \Psi^r(-z) \sigma_2,
\end{align*}
where
\begin{align*}
\sigma_2 = \begin{mat} 0 & -i \\ i & 0 \end{mat}, B(z) = \left( \zeta(z) \frac{z+1/2}{z-1/2} \right)^{\nu\sigma_3}, \zeta(z) = 2\sqrt{-t\theta(z)+t\theta(1/2)}.
\end{align*}
Consider the sectionally holomorphic matrix-valued function $\hat \Psi(z)$ defined by
\begin{align*}
\hat \Psi(z) = \begin{choices} P(z), \when |z\pm 1/2| > R,\\
\hat \Psi^r(z), \when |z - 1/2| < R,\\
\hat \Psi^l(z), \when |z + 1/2| < R. \end{choices}
\end{align*}
We use $[\hat G;\hat \Gamma]$ to denote the RHP solved by $\hat
\Psi$.  See the top panel of Figure \ref{ParametrixDeformation} for $\hat \Gamma$. In \cite{FokasPainleve} it is shown that $\Psi^r$ satisfies the RHP for $\Phi$ exactly near $z = 1/2$ and for  $\Psi^l$ near $z = -1/2$.  Notice that $\hat \Psi^r$ and $\hat \Psi^l$ are bounded near $z = \pm 1/2$.  In the special case where $\log d_1 \in \mathbb R$, $P$ remains bounded at $\pm 1/2$.  Following the analysis in \cite{FokasPainleve} we write
\begin{align*}
\Phi(z) = \chi(z) \hat \Psi(z),
\end{align*}
where $\chi \goto I$ as $\zeta \goto \infty$.

We deform the RHP for $\hat
\Psi$ to open up a small circle of radius $r$ near the origin as in Figure~\ref{RemoveParametrixCurve}.  We use $[\hat G_1;\hat \Gamma_1]$ to denote
this deformed RHP and solution $\hat \Psi_1$.  See Figure~\ref{ParametrixDeformation} for $\hat
\Gamma_1$.  Also it follows that  $\hat \Psi(z) P^{-1}(z)$ is uniformly bounded in $z$ and $\xi$. Further, $\hat \Psi_1$ has the same properties. Since $\hat \Psi_1$ is uniformly bounded in both $z$ and $\xi$ we use \eqref{explicit-inverse} to show that $\mathcal C[\hat G_1; \hat \Gamma_1]^{-1}$ has uniformly bounded norm.  We wish to use this to show the uniform boundedness of the inverse $\mathcal C[G;\Gamma]^{-1}$. To do this we extend the jump contours and jump matrices in the following way.  Set $\Gamma_e = \Gamma \cup \hat \Gamma_1$ and define
\begin{align*}
G_e(z) &= \begin{choices} G(z) \when z \in \Gamma,\\
I \otherwise,\end{choices}\\
\hat G_e(z) &= \begin{choices} \hat G_1(z) \when z \in \hat \Gamma_1,\\
I \otherwise. \end{choices}
\end{align*} 

The estimates in \cite{FokasPainleve} show that $G_e-\hat G_e \goto 0$
uniformly on $\Gamma_e$.  It follows that $\mathcal C[\hat
G_e;\Gamma_e]^{-1}$ is uniformly bounded since the extended operator
is the identity operator on $\Gamma \setminus \hat \Gamma_1$.  Theorem
\ref{approx-prop} implies that $\mathcal C[G_e;\Gamma_e]^{-1}$ is
uniformly bounded for   sufficiently large $\xi$, which implies that
$\mathcal C[G;\Gamma]^{-1}$ is uniformly bounded for $\xi$ sufficiently
large, noting that the extended operator is the identity
operator on the added contours.   We now use this construction to prove the uniform convergence of the numerical method using both direct and indirect estimates.  

\subsection{Application of Direct Estimates}

We proceed to show that the RHP for $\Psi$ satisfies the properties of a numerical parametrix.  This requires that the jump matrices have uniformly bounded Sobolev norms.  The only singularities in the jump matrices is of the form
\begin{align*}
s(z) = \left(\frac{z - 1/2}{z + 1/2}\right)^{iv}, v \in \mathbb R.
\end{align*}
After transforming to a local coordinate $k$, $z = \xi^{-1/2} k - 1/2$, we see that
\begin{align*}
S(k) = s(\xi^{-1/2} k - 1/2) = \xi^{-iv/2} \left( \frac{\xi^{-1/2} k + 1}{k} \right)^{iv}.
\end{align*}
The function $S(k)$ is smooth and has uniformly bounded derivatives provided $k$ is bounded away from $k = 0$. The deformations applied thus far guarantee that $k$ will be bounded away from $0$. To control behavior of the solution for large $k$ we look at the exponent which appears in the jump matrix
\begin{align*}
\theta(z) = \frac{2i}{3} - 4i\left(z + \half\right)^2 + \frac{8i}{3} \left(z + \half\right)^3,
\end{align*}
and define
\begin{align*}
\Theta(k) = \theta(\xi^{-1/2} k - 1/2) = \frac{2i}{3} - 4ik^2/\xi + \frac{8i}{3}k^3/\xi^{3/2}.
\end{align*}
If we assume that the contours are deformed along the {local} paths of
steepest descent, all derivatives of $e^{\xi\Theta(k)}$ are
exponentially decaying, uniformly in $\xi$.  After applying the same
procedure at $z = 1/2$ and after contour truncation, Theorem
\ref{contour-truncation} implies the RHP for $\Psi$ satisfies the hypotheses of Theorem \ref{direct-uniformity}, proving strong uniform convergence.

\subsection{Application of Indirect Estimates}

\begin{figure}[ht]
\centering
\includegraphics[width=.5\linewidth]{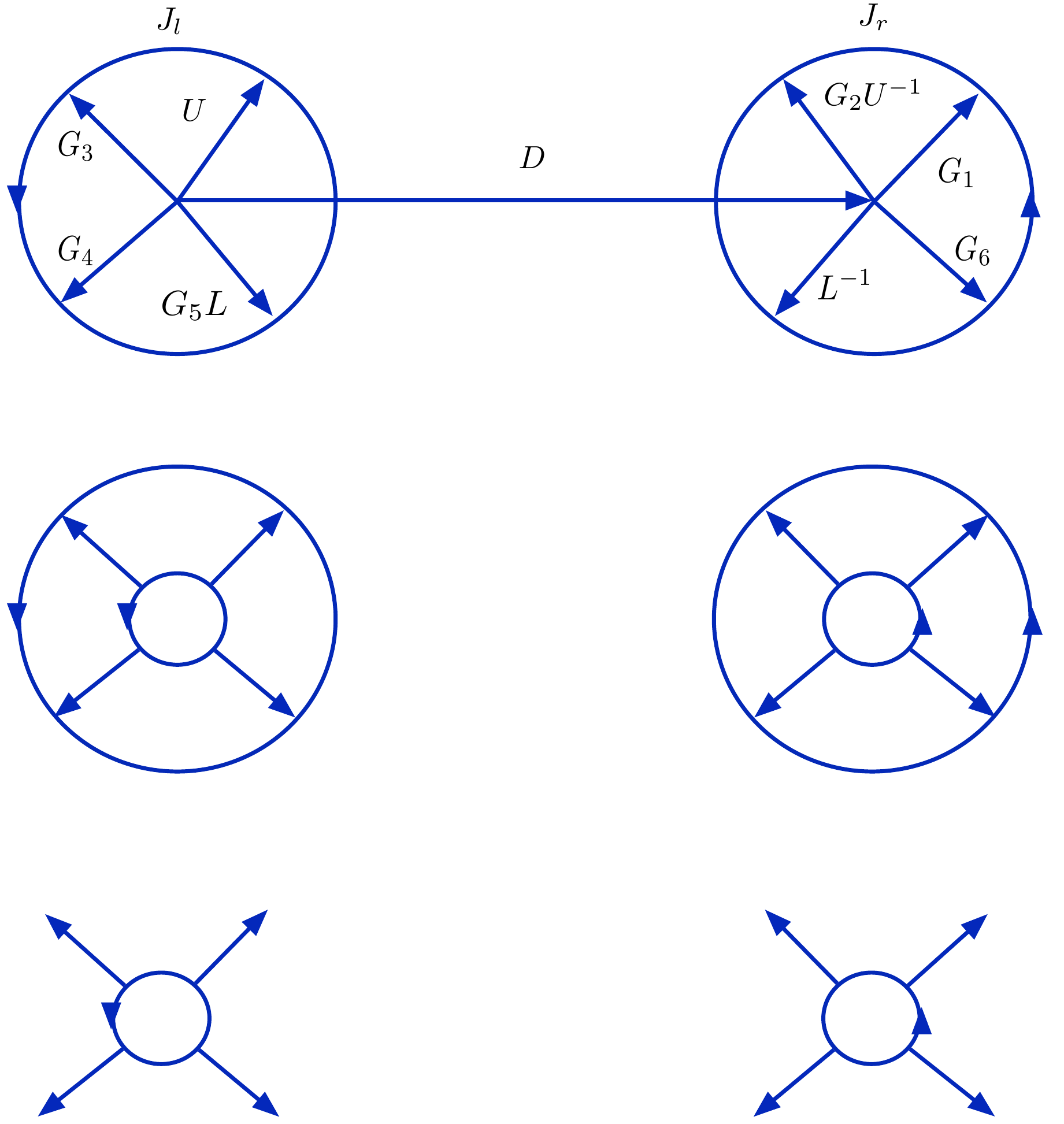}
\caption{Top:  Jump contours for the model problem with solution $\hat
  \Psi$.  Note that $J_r$ and $J_l$ are the jumps on the outside of
  the circles.  They tend uniformly to the identity as $\xi \goto
  \infty$ \cite{FokasPainleve}.  Center:  The jump contours, $\hat
  \Gamma_1$, for the function $\hat \Psi_1$. The inner circle has radius $r$ and the outer circle has radius $R$.   Bottom: Contour on which $\hat U$ is non-zero. This can be matched up with the right contour in Figure \ref{RemoveParametrixCurve}.  \label{ParametrixDeformation}}
\end{figure}

The second approach is to use the solution of the model problem to
construct an numerical parametrix.  Since we have already established
strong uniform convergence we proceed to establish a
theoretical link with the method of nonlinear steepest descent, {demonstrating that the success of  nonlinear steepest descent implies the success of the numerical method, even though the numerical method does not depend on the details of the  nonlinear steepest descent}.  We start with the RHP  $[\hat G_1;\hat \Gamma_1]$ and its solution $\hat \Psi_1$. As before, see Figure \ref{ParametrixDeformation} for $\hat
\Gamma_1$. Define $\hat u = (\hat \Psi_1)^+ - (\hat \Psi_1)^-$ which is
the solution of the associated SIE on $\hat \Gamma_1$.  The issue here
is that we cannot scale the deformed RHP in Figure
\ref{RemoveParametrixCurve} so that it is posed on the same contour as
$[G;\Gamma]$.  We need to remove the larger circle.

We define a new function $\hat U = \hat u \phi$ where $\phi$ is a
$C^\infty$ function with support in $(B(-1/2,R) \cup B(1/2,R)) \cap
\hat \Gamma_1$ such that $\phi = 1$ on $(B(1/2,r) \cup B(-1/2,r)) \cap
\hat \Gamma_1$ for $r < R$.  Let $\hat \Gamma_2$ be the support of
$\hat U$ (see bottom contour in Figure \ref{ParametrixDeformation}).   Define $\hat \Psi_2 = I + \mathcal C_{\hat
  \Gamma_2} \hat U$.  From the estimates in \cite{FokasPainleve}, it follows that
\begin{align*}
\hat \Psi_2^+ = \hat \Psi_2^- \hat G_2
\end{align*}
where $\hat G_2 - G$ tends uniformly to zero as $\xi \goto
\infty$.  We have
to establish the required smoothness of $\hat U$. We do this
explicitly from the above expression for $\hat \Psi P^{-1}$ after using the scalings $z = \xi^{-1/2} k \pm 1/2$.  The final step is to let $\xi$ be large enough so that we can truncate both $[G;\Gamma]$ and $[\hat G_2;\hat \Gamma_2]$ to the same contour.  We use Proposition \ref{local-prop} to prove that this produces a numerical parametrix.  Additionally, this shows how the local solution of RHPs can be tied to stable numerical computations of solutions.

\begin{remark}  This analysis relies heavily on the boundedness of $P$.  These arguments would fail if we were to let $P$ have unbounded singularities. In this case one approach would be to solve the RHP for $\chi$.  The jump for this RHP tends to the identity.  To prove weak uniformity for this problem one {only needs to consider} the trivial RHP with the jump being the identity matrix as a numerical parametrix.\end{remark}

\subsection{Numerical Results}

\begin{figure}[ht]
\centering
\includegraphics[width=\linewidth]{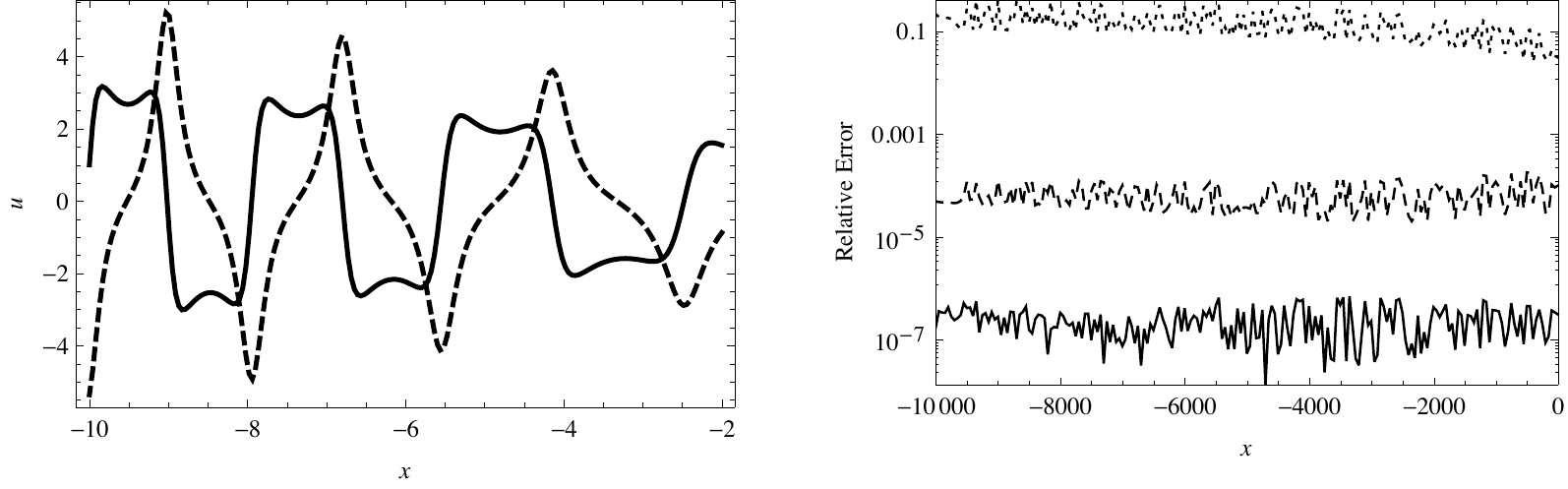}
\caption{Left:  Plot of the solution, $u$, for small $x$  (Solid: real part, Dashed: imaginary part). Right: Relative error.  Solid:  $|x|^{-1/2}|u(12,x)-u(36,x)|$, Dashed: $|u(8,x)-u(36,x)|/\sqrt{|x|}$, Dotted: $|u(4,x)-u(36,x)|/\sqrt{|x|}$.  {This plot demonstrates both uniform approximation and spectral convergence.} \label{PIINumerics}}
\end{figure}

In Figure \ref{PIINumerics} we plot the solution to \PII with $(s_1,s_2,s_3) = (1,-2,3)$ and {demonstrate numerically that the computation remains accurate in the asymptotic regime}. We use $u(n,x)$ to denote the approximate solution obtained with $n$ collocation points per contour.  Since we are using \eqref{PII-reconstruct}  we consider the estimated relative error by dividing the absolute error by $\sqrt{x}$.  We see that we retain relative error as $x$ becomes large.

\begin{remark}
	{Solutions to \PII often have poles on the real line, which correspond to the RHPs not having a solution.  In other words, $\|\mathcal C[\Gamma,\Omega]^{-1}\|$ is not uniformly bounded, which means that the theory of this paper does not apply.  However, the theorems can be adapted to the situation where  $x$ is restricted to a subdomain of the real line such that $\|\mathcal C[\Gamma,\Omega]^{-1}\|$ is uniformly bounded.  This demonstrates asymptotic stability of the numerical method for solutions with poles, provided that $x$ is bounded away from the poles,  {similar to the restriction of the asymptotic formul\ae\ in \cite{FokasPainleve}.}}  
\end{remark}

\section{Application to the modified Korteweg--de Vries Equation}\label{section:kdv}

In this section we consider the numerical solution of the modified Korteweg-de Vries equation \eqref{mKdV} (mKdV) for $x < 0$.  The RHP for mKdV is \cite{TrogdonSOKdV}
\begin{align*}
\Phi^+(z) &= \Phi^-(z) G(z), ~~z \in \mathbb R,\\
\Phi(\infty) &= I,\\
G(z) &= \begin{mat} 1-\rho(z)\rho(-z) & - \rho(-z) e^{-\theta(z)} \\ \rho(z) e^{\theta(z)} & 1 \end{mat},\\
\theta(z) &= 2izx+8iz^3t.
\end{align*}

\begin{figure}[ht]
\centering
\includegraphics[width=.5\linewidth]{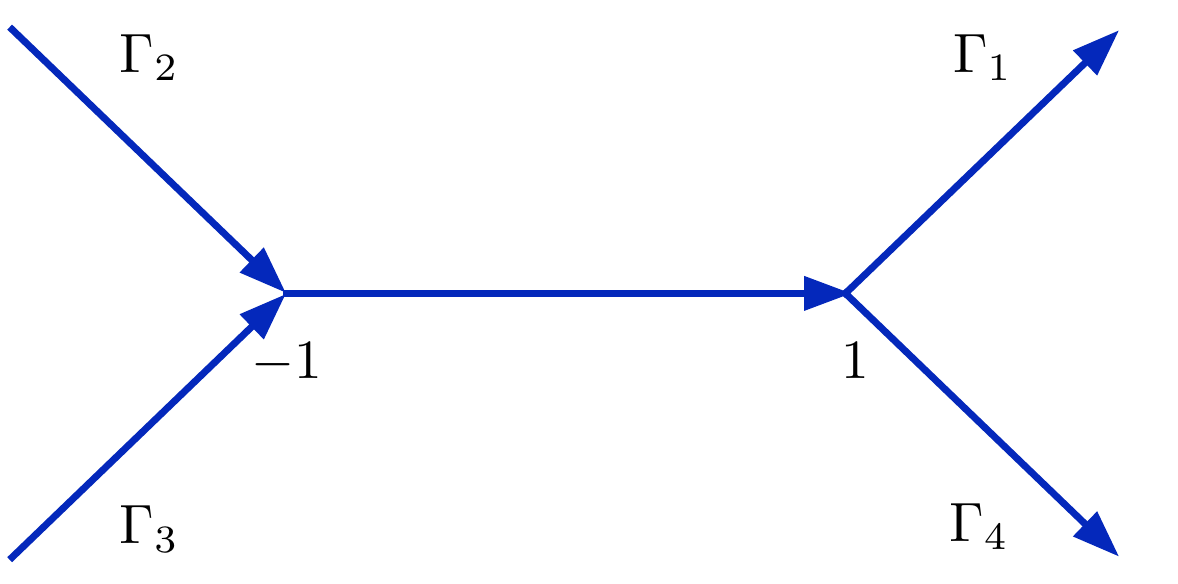}
\caption{Jump contours for the RHP for mKdV. \label{KdVPainleve}}
\end{figure}

In the cases we consider $\rho$ is analytic in a strip that contains $\mathbb R$.  If $x \ll -ct^{1/3}$ the deformation is similar to the case considered above for \PII and asymptotic stability follows by the same arguments.   We assume $x = - 12c^2 t^{1/3}$ for some positive constant $c$.  This deformation is found in \cite{TrogdonSOKdV}.  We rewrite $\theta$:
\begin{align*}
\theta(z) = -24ic^2(zt^{1/3}) + 8i (zt^{1/3})^{3}.
\end{align*}
We note that $\theta'(z_0) = 0$ for $z_0 = \pm \sqrt{-x/(12t)} = \pm c t^{-1/3}$.  We introduce a new variable $k = zt^{1/3}/c$ so that 
\begin{align*}
\theta(k ct^{-1/3}) = -24ic^3k + 8ic^3k^3 = 8ic^3(k^3-3k).
\end{align*}
For a function of $f(z)$ we use the scaling $\tilde f(k) = f(kct^{-1/3})$.  The functions $\tilde \theta$, $\tilde G$ and $\tilde \rho$ are identified similarily.  After deformation and scaling, we obtain the following RHP for $\tilde \Phi(k)$:
\begin{align*}
\tilde \Phi^+(k) = \tilde \Phi^-(k) J(k), k \in \Sigma = [-1,1] \cup \Gamma_1 \cup \Gamma_2 \cup \Gamma_3 \cup \Gamma_4,\\
J(k) = \begin{choices} \tilde G(k), \when k \in [-1,1],\\
\begin{mat} 1 & 0 \\ \tilde \rho(k)e^{\tilde \theta(k)} & 1 \end{mat}, \when k \in \Gamma_1 \cup \Gamma_2,\\
\begin{mat} 1 & -\tilde \rho(-k)e^{-\tilde \theta(k)} \\ 0 & 1 \end{mat}, \when k \in \Gamma_3 \cup \Gamma_4,
\end{choices}
\end{align*}
 where $\Gamma_i$, $i=1,2,3,4$, shown in Figure \ref{KdVPainleve}, are locally deformed along the path of steepest descent.  To reconstruct the solution to mKdV we use the {formula}
\begin{align}\label{mkdv-reconstruct}
u(x,t) = 2iz_0 \lim_{k\goto \infty} k \tilde \Phi_{12}(k).
\end{align}

\begin{remark}  We assume $\rho$ decays rapidly at $\infty$ and is analtytic in a strip that contains the real line.  This allows us to perform the initial deformation which requires modification of the contours at $\infty$. As $t$ increases, the analyticity requirements on $\rho$ are reduced; the width of the strip can be taken to be smaller if needed.  We only require that each $\Gamma_i$ lies in the domain of analticity for $\tilde \rho$.  More specifically, we assume $t$ is large enough so that when we truncate the contours for numerical purposes, they lie within the strip of analyticity for $\tilde \rho$.
\end{remark}

The parametrix derived in \cite{deift-zhou:mkdv}  is
used to show that $\mathcal C[J,\Sigma]$ has an inverse that is
uniformly bounded by using  \eqref{explicit-inverse} as was done in the previous section. We use the analyticity and decay of
$\rho$ at $\infty$ along with the fact that the contours pass along
the paths of steepest descent.

{The contour is fixed (i.e., independent of $x,t$ and $c$), and this situation is more straightforward to analyze than the previous example. }   Repeated differentiation of $J(k)$ proves that this deformation yields a uniform numerical approximation.  Furthermore, replacing $c$ by any smaller value yields the same conclusion.  This proves the uniform approximation of mKdV in the {Painlev\'e} region
\begin{align*}
 \{(x,t): t \geq \epsilon,~~ x \leq -\epsilon, x \geq -12c^2t^{1/3} \}, ~~\epsilon > 0.
\end{align*}
{where $\epsilon$ is determined by the analyticity of $\rho$.}

\subsection{Numerical Results}

\begin{figure}[ht]
\centering
\includegraphics[width=\linewidth]{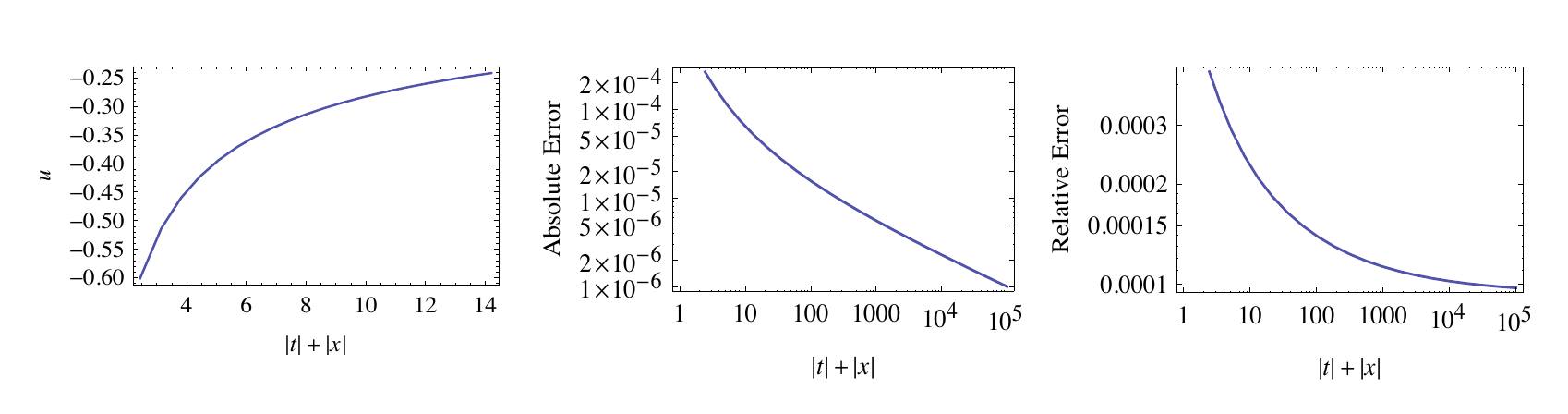}
\caption{Left:  Plot of the solution along $x = - (3t)^{1/3}$ for small time.  Center: Absolute error, $|u(5,x,t)-u(10,x,t)|$, for long time  Right:  Relative error $|u(5,x,t)-u(10,x,t)|/|u(10,x,t)$ for long time. \label{mKdVNumerics}}
\end{figure}

In Figure \ref{mKdVNumerics} we show the solution with initial data $u(x,0) = -2 e^{-x^2}$ with $c = \sqrt{9/4}$.  The reflection coefficient is obtained using the method described in \cite{TrogdonSOKdV}.  We use the notation $u(n,x,t)$ to denote the approximate solution obtained with $n$ collocation points per contour.  We see that the absolute error tends to zero rapidly. More importantly, the relative error remains small. We approximate the solution uniformly on the fixed, scaled contour.  When we compute the solution using \eqref{mkdv-reconstruct} we multiply by $z_0$ which is decaying to zero along this trajectory.  This is how the method  maintains accuracy even when comparing relative error.

\section{Acknowledgments}

We thank Bernard Deconinck for many useful conversations.  Further, support by the National Science Foundation is acknowledged through grant NSF-DMS-1008001 (TT).  Any opinions, findings, and conclusions or recommendations expressed in this material are those of the authors and do not necessarily 
reflect the views of the funding sources.

\bibliographystyle{plain}
\bibliography{NNSD}

\end{document}